\documentclass[a4paper,oneside]{amsart}
\usepackage[left=3cm,right=3cm,top=3cm,bottom=3cm]{geometry}

\usepackage{textcase}
\usepackage{amsmath}
\usepackage{amssymb}

\usepackage{xypic}
\usepackage[pagebackref]{hyperref}
\usepackage{amsthm}

\DeclareMathOperator{\supp}{supp}
\DeclareMathOperator{\rank}{rank}
\DeclareMathOperator{\Hom}{Hom}
\DeclareMathOperator{\Tor}{Tor}
\DeclareMathOperator{\Ext}{Ext}
\DeclareMathOperator{\Ker}{Ker}
\DeclareMathOperator{\Coker}{Coker}
\DeclareMathOperator{\Img}{Im}
\DeclareMathOperator{\chr}{char}

\DeclareMathOperator{\gen}{gen}
\DeclareMathOperator{\rel}{rel}

\def\id{{\mathrm{id}}}
\def\CC{\mathbb{C}}
\def\ZZ{\mathbb{Z}}

\def\QQ{\mathbb{Q}}
\def\RR{\mathbb{R}}

\def\TT{\mathbb{T}}
\def\Zg{\ZZ_{\geq 0}}
\def\Zm{\Zg^m}

\def\K{\mathcal{K}}
\def\ZK{\mathcal{Z}_\K}
\def\OZK{\Omega\ZK}
\def\DJ{\mathrm{DJ}}
\def\ODJ{\Omega\DJ}
\def\k{\mathbf{k}}
\def\b{\widetilde{b}}
\def\H{{\widetilde{H}}}
\def\kK{\k[\K]}
\def\kKc{\k\langle\K\rangle}

\def\B{\mathrm{B}}
\def\oB{\overline{\B}}

\newtheorem{thm}{Theorem}[section]
\newtheorem{lmm}[thm]{Lemma}
\newtheorem{cnj}[thm]{Conjecture}
\newtheorem{prp}[thm]{Proposition}
\newtheorem{crl}[thm]{Corollary}

\theoremstyle{definition}

\newtheorem*{rmk*}{Remark}
\newtheorem{dfn}[thm]{Definition}
\newtheorem{rmk}[thm]{Remark}
\newtheorem{algo}[thm]{Algorithm}
\newtheorem{prb}[thm]{Problem}
\numberwithin{equation}{section}

\hypersetup{
    colorlinks=true,
    linkcolor=blue,
    citecolor=blue,
    urlcolor=blue,
    pdftitle={Loop homology of moment-angle complexes in the flag case},
    pdfauthor={Fedor Vylegzhanin},
}

\title{Loop homology of moment-angle complexes in the flag case}
\author{Fedor Vylegzhanin}

\address{\parbox{\linewidth}{
Steklov Mathematical Insitute of Russian Academy of Sciences, Moscow, Russia;\\
National Research University Higher School of Economics, Russia
}}
\email{vylegf@gmail.com}


\begin{document}

\maketitle

\begin{abstract}
We develop a general homological approach to presentations of connected graded associative algebras, and apply it to the loop homology of moment-angle complexes $Z_K$ that correspond to flag simplicial complexes $K$. For an arbitrary coefficient ring, we describe generators of the Pontryagin algebra $H_*(\Omega Z_K)$ and the defining relations between them. We prove that such moment-angle complexes are coformal over $\mathbb{Q},$ give a necessary condition for rational formality, and compute their homotopy groups in terms of homotopy groups of spheres.
\end{abstract}

\section{Introduction}
For a simply connected space $X$ and a commutative ring $\k$ with unit, the \emph{Pontryagin algebra} $H_*(\Omega X;\k)$ is a connected graded associative $\k$-algebra with respect to the Pontryagin product.
We study the Pontryagin algebras of \emph{moment-angle complexes} $X=\ZK:=(D^2,S^1)^\K$ that correspond to simplicial complexes $\K$. Moment-angle complexes play an important role in toric topology \cite{ToricTopology}, and they have interesting homotopical properties and surprising connections to several topics in algebra and combinatorics \cite{bbc}. If $\K$ is a simplicial complex on the vertex set $[m]=\{1,\dots,m\},$ there is an effective action of the $m$-dimensional torus $\TT^m= (S^1)^{\times m}$ on $\ZK.$ The homotopy quotient $E\TT^m\times_{\TT^m}\ZK$ (the Borel construction) is known as the \emph{Davis--Januszkiewicz space} $\DJ(\K)$ and is homotopy equivalent to the polyhedral product $(\CC P^\infty,*)^\K,$ see \cite[Theorem 4.3.2]{ToricTopology}.

Panov and Ray \cite{pr} reduced the study of corresponding Pontryagin algebras to an algebraic problem. Applying the based loops functor to the homotopy fibration
\begin{equation}
\label{eq:zk_dj_bt_fibration}
\ZK\to\DJ(\K)\to B\TT^m,
\end{equation}
they obtained a split fibration of H-spaces $\OZK\to\ODJ(\K)\to \TT^m$
and thus an extension of cocommutative Hopf algebras
$$\k\to H_*(\OZK;\k)\to H_*(\ODJ(\K);\k)\to\Lambda[u_1,\dots,u_m]\to \k$$
over a field $\k.$ For any $\K,$ there is an isomorphism of Hopf algebras $H_*(\ODJ(\K);\k)\cong\Ext_{\kK}(\k,\k)$ \cite{pr,franz_hga} (moreover, this is true for any principal ideal domain $\k$ such that $H_*(\ODJ(\K);\k)$ is a free $\k$-module).
If $\K$ is a \emph{flag} simplicial complex, this Hopf algebra is known completely: it is isomorphic to the partially commutative algebra
$$\kK^!:=T(u_1,\dots,u_m)/(u_i^2=0,~i=1,\dots,m;~u_iu_j+u_ju_i=0,~\{i,j\}\in\K),~\deg u_i=1.$$
Generators $u_i$ are primitive and have degree $(-1,2e_i)$
with respect to the $\ZZ\times\Zm$-grading introduced in \cite{cat(zk)}.
In this case Grbi\'c, Panov, Theriault and Wu \cite{gptw} found a minimal generating set for the algebra $H_*(\OZK;\k),$ and the author calculated the number of relations in any minimal presentation (by homogeneous generators and relations) of this algebra \cite{cat(zk)}.

The last calculation relies on homological methods developed by Wall \cite{wall} and Lemaire \cite{lemaire} for connected graded associative algebras over a field. Namely, multiplicative generators of a connected $\k$-algebra $A$ correspond to additive generators of the graded $\k$-module $\Tor^A_1(\k,\k),$ and relations correspond to generators of $\Tor^A_2(\k,\k).$ In order to study the integer Pontryagin algebra $H_*(\OZK;\ZZ)$, we generalize these results to the case of \emph{arbitrary} commutative rings $\k$ with unit, and construct explicit presentations of connected $\k$-algebras using cycles in the bar construction. These results are presented in Appendix \ref{sec:appendix_presentations}. We hope that they will be useful in other contexts.\\

Let us give a general description of our approach. Suppose that we are given a connected $\k$-algebra $A$ which is a free left module over its subalgebra $S,$ $A\simeq S\otimes_\k V$. We wish to construct a presentation of $S.$ Theorem \ref{thm:presentation_from_cycles} does that, if we know a set of cycles in the bar construction $\oB(S)$, such that their images generate the $\k$-modules $H_i(\oB(S))\simeq\Tor^S_i(\k,\k),$ $i=1,2$. The following algorithm computes such cycles:
\begin{enumerate}
\item Build a free resolution $(A\otimes M,d)$ of the left $A$-module $\k.$
\item Interpret it as a free resolution $(S\otimes V\otimes M,\widehat{d})$ of the left $S$-module $\k.$ Compute the functor $\Tor^S(\k,\k)$ as the homology of the complex $(V\otimes M,\overline{d}).$ Find cycles in $(V\otimes M,\overline{d})$ such that their images generate $\Tor^S_i(\k,\k),$ $i=1,2.$
\item Construct a morphism $\varphi:(S\otimes V\otimes M,\widehat{d})\to (\B(S),d_{\B})$ of free resolutions of the left $S$-module $\k$, using the contracting homotopy of the bar resolution (see Corollary \ref{crl:contracting_homotopy_builds_resolution}). Obtain a morphism of chain complexes $\overline{\varphi}:(V\otimes M,\overline{d})\to(\oB(S),d_{\oB})$ that induces an isomorphism on the homology.
\item Applying $\overline{\varphi}$ to the cycles from (2), obtain the required cycles in $\oB(S).$
\end{enumerate}

This situation takes place if $\k\to S\to A\to V\to \k$ is an extension of connected Hopf algebras, see \cite{hopf_extension}, \cite[Proposition 4.9]{milnor_moore}. In that sense, our algorithm has similarities with the \emph{Reidemeister--Schreier algorithm} that constructs a presentation of a subgroup, given a presentation of the whole group. See \cite{hopf_subalgebras} for another approach to Hopf subalgebras in connected Hopf algebras. It is well known that extensions of Hopf algebras arise in the study of fibrations $F\to E\to B$ that have a section after looping (see Appendix \ref{sec:appendix_section} for the proof). For such ``$\Omega$-split'' fibrations, the proposed method allows to study presentations of $H_*(\Omega F;\k)$, if the algebras $H_*(\Omega E;\k)$ and $H_*(\Omega B;\k)$ are known. 

Fibrations of this kind are studied by Theriault \cite{theriault_lectures}, see also \cite[Proposition 6.1]{beben_theriault}. (However, these works deal with cases when the algebra $H_*(\Omega F;\k)$ is known better than $H_*(\Omega E;\k).$) We consider the case $F=\ZK,$ $E=\DJ(\K),$ $B=(\CC P^\infty)^m.$ The algorithm is also applicable to \emph{partial quotients} of moment-angle complexes \cite[\S 4.8]{ToricTopology} (we will consider their Pontryagin algebras in subsequent publications) and polyhedral products of the form $(PX,\Omega X)^\K$ (here we refer to the recent work \cite{cai_hopf} by Li Cai).

\subsection{Main results}
We give a presentation of the algebra $H_*(\OZK;\k)$ for a flag simplicial complex $\K$ and any ring $\k$. The presentation is explicit up to a rewriting process described in Algorithm \ref{algo:rewriting}. For $x\in H_*(\ODJ(\K);\k)$ and a subset $A=\{a_1<\dots<a_k\}\subset[m]$,  denote
$$c(A,x):=[u_{a_1},[u_{a_2},\dots[u_{a_k},x]\dots]]\in H_*(\ODJ(\K);\k).$$
This element belongs to the subalgebra $H_*(\OZK;\k)\subset H_*(\ODJ(\K);\k),$ if $x=u_i$ and $A\neq\varnothing$ (see Corollary \ref{crl:c(I,uj)_in_hozk}). For every $J\subset[m],$ denote by $\Theta(J)$ the set of all vertices $i\in J$ such that
\begin{itemize}
\item Vertices $i$ and $\max(J)$ are in different path components of the complex $\K_J;$
\item $i$ is the smallest vertex in its path component.
\end{itemize}
Denote by $\widetilde{b}_i(X;\k)$ the minimal number of elements that generate the $\k$-module $\widetilde{H}_i(X;\k).$ Clearly, $|\Theta(J)|=\widetilde{b}_0(\K_J;\k)$ for any principal ideal domain $\k$. Consider the $\widetilde{b}_0(\K_J;\k)$-element set
$$\Big\{c(J\setminus\{i\},u_i):~J\subset[m],~i\in\Theta(J)\Big\}\subset H_*(\OZK;\k).$$
We call its elements the \emph{GPTW generators} (after Grbi\'c, Panov, Theriault and Wu).
\begin{thm}
\label{thm:intro_presentation}
Let $\k$ be a commutative ring with unit and $\K$ be a flag simplicial complex without ghost vertices on vertex set $[m].$
\begin{enumerate}
\item For every $J\subset[m]$, choose a set of simplicial 1-cycles
$$\sum_{\{i<j\}\in\K_J}\lambda_{ij}^{(\alpha)}[\{i,j\}]\in C_1(\K_J;\k)$$
that generate the $\k$-module $H_1(\K_J;\k).$ Then the algebra $H_*(\OZK;\k)$ is generated by GPTW generators modulo the relations
$$
\sum_{\{i<j\}\in\K_J}\lambda_{ij}^{(\alpha)}\sum_{\begin{smallmatrix}
J\setminus\{i,j\}=A\sqcup B:\\
\max(A)>i,~\max(B)>j
\end{smallmatrix}}
\pm
\Big[\widehat{c}(A,u_i),\widehat{c}(B,u_j)\Big]=0$$
that correspond to the chosen 1-cycles. (Here $\widehat{c}(A,u_i),$ $\widehat{c}(B,u_j)$ are the elements $c(A,u_i),$ $c(B,u_j)\in H_*(\OZK;\k)$ that are arbitrarily expressed through the GPTW generators, and $[x,y]:=x\cdot y - (-1)^{|x|\cdot |y|}y\cdot x.$) In particular, $H_*(\OZK;\k)$ admits a $\ZZ\times\Zm$-homogeneous presentation by $\sum_{J\subset[m]}\b_0(\K_J;\k)$ generators and $\sum_{J\subset[m]}\b_1(\K_J;\k)$ relations.

\item If $\k$ is a principal ideal domain, then this presentation is minimal: any $\ZZ\times\Zm$-homogeneous presentation of $H_*(\OZK;\k)$ contains at least $\sum_{J\subset[m]}\b_0(\K_J;\k)$ generators and at least $\sum_{J\subset[m]}\b_1(\K_J;\k)$ relations.
\end{enumerate}
\end{thm}
This theorem follows from Theorem \ref{thm:minimal_generators} and Theorem \ref{thm:minimal_relations}, proven in Section \ref{sec:applications}. For field coefficients, these results were partially obtained in the work of Grbi\'c, Panov, Theriault, Wu (the minimal set of generators \cite[Theorem 4.3]{gptw}) and the author (number of relations and their degrees \cite[Corollary 4.5]{cat(zk)}). Sometimes the number of relations can be reduced, if we do not require them to be $\ZZ\times\Zm$-homogeneous (see Theorem \ref{thm:hozk_z_homogeneous}).

We also present new results on the homotopy of moment-angle complexes that correspond to flag complexes. Using a result of Huang \cite{huang}, we prove in Corollary \ref{crl:zk_is_coformal} that in the flag case $\ZK$ is \emph{coformal} over $\QQ$ in the sence of rational homotopy theory. Results of Berglund \cite{berglund_koszul} then give a necessary condition for such moment-angle complexes to be rationally formal (Theorem \ref{thm:zk_flag_formality_criterion}).
Finally, we improve a recent result of Stanton \cite{stanton} about the homotopy type of $\OZK$ by finding the explicit number of spheres in the product:
\begin{thm}
\label{thm:homotopy_groups_flag}
Let $\K$ be a $(d-1)$-dimensional flag simplicial complex on $[m]$ with no ghost vertices. Then there is a homotopy equivalence
\begin{equation}
\label{eq:OZK_homotopy_type}
\OZK\simeq\prod_{n\geq 3}(\Omega S^n)^{\times D_n},
\end{equation}
where the numbers $D_n\geq 0$ are determined by the identity
\begin{equation}
\label{eq:BCD_n}
-\sum_{J\subset[m]}\widetilde{\chi}(\K_J)\cdot t^{|J|}=(1+t)^{m-d}h_\K(-t)=\prod_{n\geq 3}(1-t^{n-1})^{D_n},
\end{equation}
$\widetilde{\chi}(X):=\chi(X)-1=\sum_{i\geq 0}(-1)^i\dim\H_i(X)$ is the reduced Euler characteristic and $h_\K(t):=\sum_{i=0}^d h_i(\K)\cdot t^i$ is the $h$-polynomial \cite[Definition 2.2.5]{ToricTopology} of $\K.$ In particular, for every $N\geq 1$ we have an isomorphism
\begin{equation}
\label{eq:ZK_homotopy_ranks}
\pi_N(\ZK)\simeq\bigoplus_{n=3}^N
\pi_N(S^n)^{\oplus D_n}.
\end{equation}
\end{thm}
This theorem is proved in Section \ref{sec:homotopy}. Using \eqref{eq:ZK_homotopy_ranks}, it is easy to describe the homotopy groups of corresponding Davis-Januszkiewicz spaces (using the fibration \eqref{eq:zk_dj_bt_fibration}) and partial quotients of moment-angle complexes, including quasitoric manifolds and smooth toric varieties (using similar fibrations, see \cite[Proposition 7.3.13]{ToricTopology} and \cite[\S 4]{franz}).

\subsection{Organisation of the paper} Section \ref{sec:algebra} consists of algebraic preliminaries. We highlight Corollary \ref{crl:contracting_homotopy_builds_resolution} that allows us to construct chain maps into the bar resolution. In Section \ref{sec:tortop_prelim} we recall notions from toric topology and discuss the properties of Pontryagin algebras $H_*(\ODJ(\K);\k)$ and $H_*(\OZK;\k)$. Main calculations are carried in Section \ref{sec:calculations}. In section \ref{sec:applications} we prove Theorem \ref{thm:intro_presentation} and consider an example. Section \ref{sec:homotopy} contains results about (co)formality and homotopy groups of moment-angle complexes in the flag case. In Appendix \ref{sec:appendix_presentations} we develop the homological tools for working with presentations of connected graded algebras over a commutative ring. In Appendix \ref{sec:appendix_section} we prove the following folklore fact: split fibrations of loop spaces correspond (by passing to homology) to extensions of Hopf algebras. Appendix \ref{sec:appendix_identities} contains commutator identities that are used in Section \ref{sec:calculations}.

\subsection{Acknowledgements}
The author thanks his advisor Taras Panov for guidance and attention to this work, Matthias Franz for pointing out \cite[Proposition 6.5]{franz_hga}, Lewis Stanton for proving Lemma \ref{lmm:eliminating_spheres} which greatly simplifies the statement of Theorem \ref{thm:homotopy_groups_flag}, the anonymous referee for important suggestions and corrections, and G. Chernykh, V. Gorchakov, D. Piontkovski, T. Rahmatullaev, and A. Saigak for useful comments and conversations.

\subsection{Funding}
This work is supported by the Russian Science Foundation under grant no. 23-11-00143, \href{http://rscf.ru/en/project/23-11-00143/}{http://rscf.ru/en/project/23-11-00143/}.
The author is a laureate of the all-Russia mathematical August Moebius contest of graduate and undergraduate student papers and thanks the jury and the board for the high praise of his work.

\section{Preliminaries: algebra}
\label{sec:algebra}
\subsection{Connected graded algebras}
Fix a commutative associative ring $\k$ with unit. We consider associative $\k$-algebras with unit that are graded by a commutative monoid $G$ (usually $G=\ZZ$ or $\ZZ^k\times \Zm,$ $k=0,1,2.$) Left $A$-modules are also $G$-graded. Elements of $\Zm$ are denoted by $\alpha=(\alpha_1,\dots,\alpha_m)=\sum_{j=1}^m\alpha_je_j,$ $\alpha_j\geq 0.$ Subsets $J\subset[m]$ are identified with elements $\sum_{j\in J}e_j\in\Zm.$ Denote also
$$|\alpha|:=\alpha_1+\dots+\alpha_m,\quad \supp\alpha:=\{i\in[m]:~\alpha_i>0\}.$$ Every $\ZZ^k\times\Zm$-graded algebra $A$ is considered as $\ZZ$-graded with respect to the total grading $A_n:=\bigoplus_{n=i_1+\dots+i_k+|\alpha|}A_{i_1,\dots,i_k,\alpha}.$

Graded algebra $A$ is \emph{connected} if $A_{<0}=0$ and $A_0=\k\cdot 1.$ We have the canonical augmentation
$\varepsilon:A\to A_0=\k$ and the augmentation ideal $I(A):=\Ker\varepsilon.$ Examples of connected $\k$-algebras:
\begin{itemize}
\item exterior algebra $\Lambda[m]:=\Lambda[u_1,\dots,u_m],$ $\deg(u_i)=(-1,2e_i)\in\ZZ\times\Zm$ with the basis $\{u_I:=u_{i_1}\wedge\dots\wedge u_{i_k},~I=\{i_1<\dots<i_k\}\};$
\item polynomial algebra $\k[m]:=\k[v_1,\dots,v_m],$ $\deg(v_i)=(0,2e_i)$ with the basis $\{v^\alpha:=\prod_{i=1}^mv_i^{\alpha_i},~\alpha\in\Zm\};$
\item tensor algebra $T(x_1,\dots,x_N),$ where $x_i$ are homogeneous elements of arbitrary positive degrees.
\end{itemize}

For a homogeneous element $a$, denote $\overline{a}:=(-1)^{1+\deg(a)}\cdot a.$ Clearly, $\overline{a\cdot b}=-\overline{a}\cdot\overline{b}$ and $\overline{\overline{a}}=a.$ 

Let $A$ be a $G$-graded algebra. Complexes of  $A$-modules $(M,d)$ are considered as $\ZZ\times G$-graded modules with a differential of degree $(-1,0).$ We use the Koszul sign rule with respect to the total grading: $d(a\cdot m)=(-1)^{\deg(a)}a\cdot d(m)=-\overline{a}\cdot d(m).$ Several formulas from \cite{cat(zk)} do not follow this rule and are corrected in this paper.

\subsection{Bar resolution and bar construction}
Let $A$ be a connected  $\k$-algebra and $\varepsilon:A\to\k$ be the augmentation. The resulting left $A$-module $\k$ has the \emph{bar resolution}
$$\dots\to \B_2(A)\to\B_1(A)\to\B_0(A)\to\k\to 0,$$
where $\B_n(A):=A\otimes I(A)^{\otimes n}.$ An element of the form $a\otimes a_1\dots\otimes a_n\in\B_n(A)$ has bidegree $(n,\deg(a)+\sum_{i=1}^n\deg(a_i))$ and is traditionally written as $a[a_1|\dots|a_n].$ The differential $d_\B$ has bidegree $(-1,0)$ and is given by the formula
$$-d_\B(a[a_1|\dots|a_n]):= \overline{a}\cdot a_1[a_2|\dots|a_n]+\sum_{i=1}^{n-1}\overline{a}[\overline{a}_1|\dots|\overline{a}_{i-1}|\overline{a}_i\cdot a_{i+1}|a_{i+2}|\dots|a_n].$$ 
Consider also the contracting homotopy $s_n:\B_n(A)\to\B_{n+1}(A),$
\begin{equation}
\label{eq:bar_contraction}
s_n(a[a_1|\dots|a_n]):=\begin{cases}
0,&a\in A_0\simeq \k;\\
[a|a_1|\dots|a_n],&\deg(a)>0;
\end{cases}
\quad
s_{-1}:\k\to\B_0(A),\quad 1\mapsto 1[].
\end{equation}
It is easy to show that $s\circ d_\B+d_\B\circ s=\id,$ $d_\B^2=0.$ Hence $(\B(A),d_{\B})$ is a free resolution of the left $A$-module $\k$, assuming that $A$ is a free $\k$-module. In this case, we obtain
$$\Tor^A_n(\k,\k)\cong H_n\Big[\oB(A),d_{\oB}\Big],$$
where $\oB(A):=\k\otimes_A\B(A)$ is the \emph{bar construction} of $A.$ We have
$$\oB_n(A) = I(A)^{\otimes n},\quad \deg([a_1|\dots|a_n]) = (n,\deg(a_1)+\dots+\deg(a_n)),~\deg d_{\oB}=(-1,0),$$
\begin{equation}
\label{eq:bar_construction}
d_{\oB}([a_1|\dots|a_n])=\sum_{i=1}^{n-1}[\overline{a}_1|\dots|\overline{a}_{i-1}|\overline{a}_i\cdot a_{i+1}|a_{i+2}|\dots|a_n]\in\oB_{n-1}(A).
\end{equation}
In particular, $d_{\oB}([x|y])=[\overline{x}\cdot y]$ and $d_{\oB}([x|y|z])=[\overline{x}\cdot y|z]+[\overline{x}|\overline{y}\cdot z].$

\subsection{Chain maps into resolutions with a contracting homotopy}
Any map of modules can be extended to a map of their free resolutions. Moreover, this extension can be described in terms of the contracting homotopy for the latter resolution. This recursive construction seems to be known to specialists: its generalisations and applications are discussed in \cite{explicit_acyclic}. The author thanks Georgy Chernykh for the reference. 

\begin{lmm}
\label{lmm:building_resolution}
Let $A$ be an associative $\k$-algebra. Suppose that the commutative diagram of left $A$-modules and their homomorphisms
$$\xymatrix{
C_n\ar[r]^-{\widehat{d}_n} & C_{n-1}\ar[r]^-{\widehat{d}_{n-1}}\ar[d]^-{\varphi_{n-1}} & C_{n-2}\ar[d]^-{\varphi_{n-2}}\\
B_n\ar[r]_-{d_n} & B_{n-1}\ar[r]_{d_{n-1}} & B_{n-2}
}$$
satisfy the conditions:
\begin{enumerate}
\item $C_n$ is a free  $A$-module with a basis $\{e_i\};$
\item $\widehat{d}_{n-1}\circ\widehat{d}_n=0;$
\item there are $\k$-linear maps $s_{n-1}:B_{n-1}\to B_n$ and $s_{n-2}:B_{n-2}\to B_{n-1}$ such that $d_n\circ s_{n-1}+s_{n-2}\circ d_{n-1}=\id_{B_{n-1}}.$
\end{enumerate}
Define an $A$-linear map $\varphi_n:C_n\to B_n$ on the basis by the formula
$$\varphi_n(e_i):= s_{n-1}(\varphi_{n-1}(\widehat{d}_n(e_i)))\in B_n.$$
Then $d_n\circ \varphi_n=\varphi_{n-1}\circ\widehat{d}_n.$
\end{lmm}
\begin{proof}
Since $d_n\circ \varphi_n$ and $\varphi_{n-1}\circ\widehat{d}_n$ are maps of $A$-modules, it is sufficient to show that they agree on the basis of $C_n.$ By definition,
$$d_n(\varphi_n(e_i))=(d_n\circ s_{n-1}\circ \varphi_{n-1}\circ \widehat{d}_n) (e_i).$$
Condition (3) gives $d\circ s\circ\varphi\circ\widehat{d}=\varphi\circ\widehat{d}-s\circ d\circ \varphi\circ\widehat{d}.$ From the commutativity of the diagram and condition (2) we obtain $s\circ d\circ\varphi\circ\widehat{d}=s\circ\varphi\circ\widehat{d}\circ\widehat{d}=0.$
Hence
$d_n(\varphi_n(e_i))= \varphi_{n-1}(\widehat{d}_n(e_i))-0.$
\end{proof}

\begin{crl}
\label{crl:contracting_homotopy_builds_resolution}
Let $A$ be a connected $\k$-algebra, $(A\otimes V_\bullet,\widehat{d}_\bullet)$ be a free resolution of the left $A$-module $\k.$ Let $\overline{\varphi}_0:V_0\to \k$ be a map of $\k$-modules such that the diagram 
$$\xymatrix{
A\otimes V_1\ar[r]^-{\widehat{d}_1}
&
A\otimes V_0\ar[r]^-{\widehat{d}_0}
\ar[d]_-{\id\otimes\overline{\varphi}_0}
&
\k
\ar@{=}[d]^-{\id=:\varphi_{-1}}
\\
\B_1(A)\ar[r]^-{d_{\B,1}}
&
A
\ar[r]^-{\varepsilon}
&
\k
}$$
commutes. Choose bases $\{e_i^{(n)}\}$ of $\k$-modules $V_n,$ and define $A$-linear maps $\varphi_n:A\otimes V_n\to \B_n(A)$ recursively as
$$
\varphi_0:=\id_A\otimes\overline{\varphi}_0,\quad
\varphi_n(a\otimes e_i^{(n)}):= a\cdot s_{n-1}(\varphi_{n-1}(\widehat{d}_n(e_i^{(n)}))),$$
where $s_{n-1}:\B_{n-1}(A)\to \B_n(A)$ is the contracting homotopy \eqref{eq:bar_contraction}.

Then $\varphi_\bullet:(A\otimes V_\bullet,\widehat{d})\to(\B_\bullet(A),d_{\B})$ is a chain map.
\end{crl}
\begin{proof}
Induction on $n.$ For $n=0$ the identity $d_{\B,n}\circ \varphi_n=\varphi_{n-1}\circ\widehat{d}_n$ holds, since the diagram commutes. The inductive step from $n-1$ to $n$ is supplied by Lemma \ref{lmm:building_resolution}.
\end{proof}
\subsection{Hopf algebra extensions and loop homology}
\label{subsec:hopf_extension}
If $A$ is a Hopf algebra over $\k$, we denote the comultiplication by $\Delta:A\to A\otimes A$ and the (co)unit maps by $\eta_A:\k\to A,$ $\varepsilon_A:A\to\k.$ Graded $\k$-Hopf algebra $A$ is \emph{connected} if $A_{<0}=0,$ $A_0=\k\cdot 1.$ The counit is then the standard augmentation $\varepsilon:A\to A_0\simeq\k.$

\begin{dfn}
Let $\iota:A\to C,$ $\pi:C\to B$ be morphisms of $\k$-Hopf algebras. They form an \emph{extension of Hopf algebras}, or a \emph{short exact sequence of Hopf algebras}
$$\k\to A\overset\iota\longrightarrow C\overset\pi\longrightarrow B\to\k,$$
if the following conditions are satisfied:
\begin{enumerate}
\item $\iota$ is injective;
\item $\pi$ is surjective;
\item $\pi\circ\iota=\varepsilon;$
\item $\Ker\pi=I(A)\cdot C;$
\item $\Img\iota=\{x\in C:~((\id_C\otimes\pi)\circ\Delta)(x)=x\otimes 1\}.$
\end{enumerate}
\end{dfn}
See \cite[Definition 1.2.0, Proposition 1.2.3]{hopf_extension} for an equivalent and more ``symmetrical'' definition.
Extensions of connected Hopf algebras were studied implicitly in \cite[\S 4]{milnor_moore}.
\begin{prp}[see {\cite[Proposition 4.9]{milnor_moore}}]
\label{prp:hopf_extension_criterion}
Let $\iota:A\to C,$ $\pi:C\to B$ be maps of connected $\k$-Hopf algebras. Suppose that a map $\Phi:A\otimes B\to C$ is an isomorphism of left $A$-modules and right $C$-comodules, and suppose that
$$\iota = \Phi\circ(\id_A\otimes\eta_B),\quad \pi\circ\Phi=\varepsilon_A\otimes\id_B.$$
Then $\k\to A\overset\iota\longrightarrow C\overset\pi\longrightarrow B\to\k$ is an extension of Hopf algebras. Conversely, for every Hopf algebra extension there is a map $\Phi$ with described properties.\qed
\end{prp}

Our main example of Hopf algebras are Pontryagin algebras (loop homology) of connected topological spaces. Let $\k$ be a commutative ring, $Y$ be a topological space such that $H_*(Y;\k)$ is a free $\k$-module. Then $H_*(Y;\k)$ is supplied with the cocommutative \emph{cup coproduct} which is dual to the cup product on $H^*(Y;\k)$: it is the composition
$$\xymatrix{
H_*(Y;\k)\ar[r]^-{\Delta_*}
&
H_*(Y\times Y;\k)
\ar[r]^-{AW_*}_-\simeq
&
H_*(C_*(Y;\k)\otimes C_*(Y;\k))
&
H_*(Y;\k)\otimes H_*(Y;\k)
\ar[l]^-\simeq_-{\kappa},
}$$
where $AW$ is the Alexander-Whitney map and $\kappa$ is the K\"unneth isomorphism. If $Y$ is also an H-space, the cup coproduct respects the Pontryagin product
$$\xymatrix{
m:H_*(Y;\k)\otimes H_*(Y;\k)
\ar[r]^-\times
&
H_*(Y\times Y;\k)
\ar[r]^-{\mu_*}
&
H_*(Y;\k)
}$$
and hence $H_*(Y;\k)$ is a cocommutative Hopf algebra. In particular, $H_*(\Omega X;\k)$ is a connected cocommutative $\k$-Hopf algebra whenever $X$ is a simply connected space such that $H_*(\Omega X;\k)$ \emph{is free over $\k$}
\cite[8.9]{milnor_moore}. Otherwise $\kappa$ fails to be an isomorphism, hence the coproduct is not defined and $H_*(\Omega X;\k)$ is merely a connected associative $\k$-algebra with unit.

In Appendix \ref{sec:appendix_section} we describe a situation when a fibration $F\to E\to B$ of simply connected spaces gives rise to an extension $\k\to H_*(\Omega F;\k)\to H_*(\Omega E;\k)\to H_*(\Omega B;\k)\to\k$ of connected Hopf algebras.

\section{Preliminaries: toric topology}
\label{sec:tortop_prelim}
\subsection{Simplicial complexes and polyhedral products}
\emph{Simplicial complex} $\K$ on the vertex set $W$ is a non-empty family of subsets $I\subset W$ that is closed under taking subsets. Elements $I\in\K$ are called \emph{faces.} We suppose that $\K$ \emph{has no ghost verties}, i.e. $\{i\}\in\K$ for all $i\in W.$ Usually $W\subset[m]:=\{1,\dots,m\}.$ Sometimes by properties of a complex $\K$ we mean properties of its \emph{geometrical realisation}, of the topological space $|\K|:=\bigcup_{I\in\K}\Delta_I\subset\Delta_W.$

For every $J\subset W$, a simplicial complex $\K_J:=\{I\in\K:~I\subset J\}$ on the vertex set $J$ (a \emph{full subcomplex} of $\K$) is defined.

Throughout the text, we write $I\setminus i:=I\setminus\{i\}$ for $i\in I$ and $I\sqcup i:=I\sqcup\{i\}$ for $i\in W\setminus I.$ Subset $I\subset W$ is a \emph{missing face} of $\K$ if $I\notin\K,$ but $I\setminus i \in\K$ for all $i\in I.$ Simplicial complex $\K$ is \emph{flag} if all its missing faces consist of two elements. 

For every complex $\K$ on vertex set $[m]$, the $\Zm$-graded \emph{Stanley--Reisner ring}
$$\kK:=\k[v_1,\dots,v_m]/\left(\prod_{i\in I}v_i=0,~I\notin\K\right),\quad \deg v_i:=2e_i\in\Zm$$
is defined. It has a homogeneous basis $\{v^\alpha:=\prod_{i=1}^m v_i^{\alpha_i}\mid\supp\alpha\in\K\}$ as a $\k$-module.
The dual $\k$-module $\kKc$ is called the \emph{Stanley--Reisner coalgebra}. It has an additive basis $\{\chi_\alpha\mid\supp\alpha\in\K\},$ $\deg\chi_\alpha=2\alpha,$ and commutative associative comultiplication
$\Delta\chi_\alpha:=\sum_{\alpha=\beta+\gamma}\chi_\beta\otimes\chi_\gamma.$

Now let $\K$ be a simplicial complex on $[m]$ and $(\underline{X},\underline{A}):=((X_1,A_1),\dots,(X_m,A_m))$ be a sequence of pairs of topological spaces. Their \emph{polyhedral product} $(\underline{X},\underline{A})^\K$ is the union
$$(\underline{X},\underline{A})^\K:=
\bigcup_{I\in\K}(\underline{X},\underline{A})^I\subset X^m,\quad (\underline{X},\underline{A})^I=Y_1\times\dots\times Y_m,\quad Y_j:=\begin{cases}
X_j,&j\in I;\\
A_j,&j\notin I.
\end{cases}$$
The addition of a ghost vertex $v$ to $\K$ replaces the space $(\underline{X},\underline{A})^\K$ with $(\underline{X},\underline{A})^\K\times A_v.$ Hence in many cases it is sufficient to consider only complexes without ghost vertices. 

Denote $(X,A)^\K:=(\underline{X},\underline{A})^\K$ if $X_i=X,$ $A_i=A$ for all $i\in[m].$ We consider two special cases of this construction: \emph{moment-angle complexes} $\ZK:=(D^2,S^1)^\K$ and \emph{Davis-Januszkiewicz spaces} $\DJ(\K):=(\CC P^\infty,*)^\K.$ It is well known that $H^*(\DJ(\K);\k)\cong\kK$ and $H^*(\ZK;\k)\cong\Tor^{\k[m]}(\kK,\k)$ as graded rings. Moreover,
$$H^n(\ZK;\k)=\bigoplus_{n=-i+2|J|}H^{-i,2J}(\ZK;\k),\quad H^{-i,2J}(\ZK;\k)\cong\H^{|J|-i-1}(\K_J;\k),$$
and the product has a geometric description in terms of maps $\K_{I\sqcup J}\hookrightarrow \K_I\ast\K_J,$ see \cite[Theorem 4.5.8]{ToricTopology}.
\subsection{Loop homology as Hopf algebras}
\begin{prp}[{\cite[Theorem 4.3.2, \S 8.4]{ToricTopology}}]
\label{prp:zk_dj_bt_fibration}
There is a homotopy fibration $\ZK\to\DJ(\K)\overset{i}\longrightarrow (\CC P^\infty)^m$ of simply connected spaces, where $i$ is the standard inclusion. The map $\Omega i$ admits a homotopy section $\sigma:\TT^m\to\Omega\DJ(\K)$ that corresponds to the choice of generators in $\pi_2(\DJ(\K))\cong\ZZ^m$ and gives rise to a homotopy equivalence $\ODJ(\K)\simeq\OZK\times\TT^m.$\qed
\end{prp}
The following description of $H_*(\ODJ(\K);\k)$ was first given in \cite[(8.4)]{pr} for $\k=\QQ,$ but the argument is easily generalised to the arbitrary coefficient ring. The main ingredients are integral \emph{formality} of $\DJ(\K)$ \cite{notbohm_ray}, Adams' cobar construction (see \cite{adams_cobar}) and a result of Fr\"oberg \cite{froeberg}.
\begin{thm}[{\cite[Theorem 1.1]{cat(zk)}}]
\label{thm:hodj_description}
For any simplicial complex $\K$ with no ghost vertices and any commutative ring $\k,$ we have an isomorphism $H_*(\ODJ(\K);\k)\cong\Ext_{\kK}(\k,\k)$  of graded $\k$-algebras (with respect to Pontryagin product and to Yoneda product). More precisely,
$$H_n(\ODJ(\K);\k)\cong\bigoplus_{-i+2|\alpha|=n}\Ext_{\kK}^i(\k,\k)_{2\alpha}.$$

This isomorphism defines the $\ZZ\times\Zm$-grading on $H_*(\ODJ(\K);\k).$ The ``diagonal'' subalgebra $D=\bigoplus_{\alpha\in\Zm}H_{-|\alpha|,2\alpha}(\ODJ(\K);\k)\subset H_*(\ODJ(\K);\k)$ is isomorphic to the algebra
$$\kK^!:=T(u_1,\dots,u_m)/(u_i^2=0,~i=1,\dots,m;~u_iu_j+u_ju_i=0,~\{i,j\}\in\K),\quad\deg u_i=(-1,2e_i).$$
For flag $\K,$ the algebra $H_*(\ODJ(\K);\k)$ coincides with $D,$ and we have $H_*(\ODJ(\K);\k)\cong\kK^!.\qed$
\end{thm}

If $H_*(\Omega Y;\k)$ is a free $\k$-module, the cup coproduct is compatible with the Pontryagin product, hence this associative algebra is a cocommutative $\k$-Hopf algebra.
Similarly, if $A$ is a commutative graded $\k$-algebra such that $\Ext_A(\k,\k)$ is a free $\k$-module, then the shuffle product on the bar construction (see \cite[Theorem X.12.2]{maclane}) induces a commutative coproduct on $\Ext_A(\k,\k)$ that is compatible with the Yoneda product. In our case, these coproduct coincide. This follows from a stronger formality result for Davis-Januszkiewicz spaces, the \emph{hga formality} \cite[Theorem 1.3]{franz_hga}.

\begin{prp}[{\cite[Proposition 6.5]{franz_hga}}]
\label{prp:hodj_is_ext}
Let $\K$ be a simplicial complex with no ghost vertices, and let $\k$ be a principal ideal domain such that $H_*(\ODJ(\K);\k)$ is a free $\k$-module. Then $H_*(\ODJ(\K);\k)\cong\Ext_{\kK}(\k,\k)$ as Hopf algebras.
\end{prp}
\begin{proof}[Outline of the proof]
Let $A$ be a dga algebra. The \emph{homotopy Gerstenhaber algebra} (hga) structure on A is a multiplication on its bar construction $\oB(A)$ such that $\oB(A)$ becomes a dga bialgebra \cite[\S 4]{franz_hga}. This structure arises naturally if $A$ is commutative (then the multiplication is the shuffle product) or if $A=C^*(X;\k)$ is the dga algebra of cochains of a 1-reduced simplicial set (then the multiplication was essentially constructed by Baues \cite[\S 2]{baues}). Then $H^*(\Omega X;\k)\cong H^*\Big[\oB(C^*(X;\k))\Big]$ as bialgebras. By a result of Franz \cite[Theorem 1.3]{franz_hga}, hga algebras $C^*(\DJ(\K);\k)$ and $\kK$ are quasi-isomorphic. The functor $\oB$ preserves quasi-isomorphisms, so $H^*(\ODJ(\K);\k)\cong H^*(\oB(\kK);\k)\cong \Tor^{\kK}(\k,\k)$ as bialgebras. Since the Hopf algebra structure on a bialgebra is unique, it is an isomorphism of Hopf algebras. The statement for $H_*(\ODJ(\K);\k)$ follows by dualisation. 
\end{proof}

\begin{rmk}
The algebra $H_*(\ODJ(\K);\k)$ is not always a free $\k$-module. For example, let $\K$ be a minimal triangulation of $\RR P^2$. Then $\ZK$ is a wedge of $\Sigma^7\RR P^2$ and spheres \cite[Example 3.3]{gptw}. We have $\ODJ(\K)\simeq\OZK\times\TT^m,$ hence $\Omega\Sigma^7\RR P^2$ is a retract of $\ODJ(\K).$ It follows that $H_*(\ODJ(\K);\ZZ)$ has 2-torsion.
\end{rmk}

Recall that an element $x$ is called \emph{primitive} if $\Delta x=x\otimes 1+1\otimes x,$ and a Hopf algebra is \emph{primitively generated} if it is multiplicatively generated by its primitive elements. 
\begin{cnj}
\label{cnj:hodj_is_primgen}
The Hopf algebra $H_*(\ODJ(\K);\k)$ is primitively generated for every simplicial complex $\K$ and every ring $\k$ such that $H_*(\ODJ(\K);\k)$ is a free $\k$-module.
\end{cnj}
By deep results of Andr\'e and Sj\"odin (see \cite[Theorem 10.2.1(5)]{avramov}), for every field $\k$ the Hopf algebra $\Ext_A(\k,\k)$ is the universal enveloping of a Lie algebra (of a 2-restricted Lie algebra, if $\chr\k=2$). In particular, this Hopf algebra is primitively generated. (This also follows from results of Browder \cite{browder}, see \cite[Theorem 10.4]{neisendorfer}.) Hence Conjecture \ref{cnj:hodj_is_primgen} holds if $\k$ is a field.
\begin{rmk}
The Hopf algebra $H_*(\Omega X;\k)$ is not always primitively generated, even if $X$ is a suspension. For example, one can take $X=\Sigma\CC P^2,$ $\k=\ZZ$ or $\ZZ/2$ (see \cite[\S 4.2]{buchstaber_grbic}). On the other hand, \cite[Theorem B]{halperin} implies that $H_*(\Omega \Sigma \CC P^d;\ZZ/p)$ is primitively generated for $p>d.$
\end{rmk}

Now we describe the connection between the loop homology of Davis-Januszkiewicz spaces and of moment-angle complexes in the form of a Hopf algebra extension.
\begin{prp}
\label{prp:ozk_extension_statement}
Let $\K$ be a simplicial complex on $[m]$ and $\k$ be a commutative ring with unit, such that $H_*(\OZK;\k)$ is a free $\k$-module. Then
$$\k\to
H_*(\OZK;\k)
\overset{\iota}\longrightarrow
H_*(\ODJ(\K);\k)\overset{p}\longrightarrow
\Lambda[u_1,\dots,u_m]
\to\k$$
is an extension of connected $\ZZ\times\Zm$-graded $\k$-Hopf algebras. The projection $p$ maps $u_i$ to $u_i.$ Its $\k$-linear section $\sigma_*:\Lambda[u_1,\dots,u_m]\to H_*(\ODJ(\K);\k)$ is given by the formula $$\sigma_*(u_I)=\widehat{u}_I:=u_{i_1}\cdot\dotso\cdot u_{i_k},\quad I=\{i_1<\dots<i_k\}.$$ Therefore, the formula $\Phi(a\otimes u_I):=\iota(a)\cdot\widehat{u}_I$ defines an isomorphism of left $H_*(\OZK;\k)$-modules and right $\Lambda[u_1,\dots,u_m]$-comodules $\Phi:H_*(\OZK;\k)\otimes\Lambda[u_1,\dots,u_m]\to H_*(\ODJ(\K);\k).$ 
\end{prp}
\begin{proof}
By Theorem \ref{thm:section_hopf_extension}, the fibration from Proposition \ref{prp:zk_dj_bt_fibration} gives rise to the required Hopf algebra extension. The formula for $p$ follows from functoriality, since the map $\DJ(\K)\hookrightarrow\DJ(\Delta_{[m]})\cong(\CC P^\infty)^m$ is induced by the inclusion $\K\hookrightarrow\Delta_{[m]}.$ The formula for $\sigma_*$ follows from the description of the homotopy section $\sigma:\TT^m\simeq \Omega B\TT^m=(\Omega\CC P^\infty)^{\times m}\to\ODJ(\K)$ as a concatenation of loops, $(\gamma_1,\dots,\gamma_m)\mapsto\gamma_1\cdot\dotso\cdot\gamma_m.$ The maps $p$ and $\sigma_*$ respect the $\ZZ\times\Zm$-grading, hence the multigrading on $H_*(\OZK;\k)$ is well defined.
\end{proof}
Since $\iota$ is injective, we identify elements of $H_*(\OZK;\k)$ with their images in $H_*(\ODJ(\K);\k).$ Let us describe some of these elements. Recall that we denote $[a,b]:=ab+(-1)^{\deg(a)\deg(b)+1}\,ba$ and $c(I,x):=[u_{i_1},[u_{i_2},\dots,[u_{i_k},x]\dots]]\in H_*(\ODJ(\K);\k)$ for $I=\{i_1<\dots<i_k\}$ and $x\in H_*(\ODJ(\K);\k).$ 
 In particular, $c(\varnothing,x):=x$ and $c(\{i\},u_j)=[u_i,u_j]=u_iu_j+u_ju_i.$
\begin{crl}
\label{crl:prim_in_hozk}
Let $x\in H_*(\ODJ(\K);\k)$ be a primitive element such that $p(x)=0.$ Then $x\in H_*(\OZK;\k).$
\end{crl}
\begin{proof}
Follows from Corollary \ref{crl:primitive_in_kernel} applied to the Hopf algebra extension from Proposition \ref{prp:ozk_extension_statement}.
\end{proof}
\begin{crl}
\label{crl:comm_in_hozk}
Let $x\in H_*(\ODJ(\K);\k)$ be a primitive element and $I\subset[m],$ $I\neq\varnothing.$ Then $c(I,x)\in H_*(\OZK;\k).$
\end{crl}
\begin{proof}
Elements $u_1,\dots,u_m\in H_*(\ODJ(\K);\k)$ are primitive for dimension reasons. Primitive elements form a Lie algebra, hence $c(I,x)\in H_*(\ODJ(\K);\k)$ is primitive.
We have $p(c(I,x))=c(I,p(x))=0,$ since it is a commutator in the commutative algebra $\Lambda[m].$ Then $c(I,x)\in H_*(\OZK;\k)$ by Corollary \ref{crl:prim_in_hozk}.
\end{proof}
\begin{crl}
\label{crl:c(I,uj)_in_hozk}
Let $j\in[m]$ and $I\subset[m],$ $I\neq\varnothing.$ Then $c(I,u_j)\in H_*(\OZK;\k).$\qed
\end{crl}

\subsection{The flag case}
Let $\K$ be a flag complex with no ghost vertices. By Theorem \ref{thm:hodj_description}, $H_*(\ODJ(\K);\k)\cong\kK^!$ is a free $\k$-module, hence the Hopf algebra structure on $H_*(\ODJ(\K);\k)$ is well defined. Moreover, the connected $\k$-algebra $\kK^!$ is generated by elements of degree $1.$ These conditions determine the Hopf algebra structure on $\kK^!$ uniquely: the elements $u_1,\dots,u_m$ are primitive. Therefore, in the flag case Conjecture \ref{cnj:hodj_is_primgen} is true for any $\k$.

The following important result was recently obtained by Stanton.
\begin{thm}[{\cite[Corollary 1.5]{stanton}}]
\label{thm:stanton_result}
Let $\K$ be a flag simplicial complex or a skeleton of a flag complex. Then $\OZK$ is homotopy equivalent to a finite type product of spaces of the form $S^1,$ $S^3,$ $S^7$ and $\Omega S^n$ for $n\geq 2,$ $n\neq 2,4,8.$\qed
\end{thm}

This gives a short proof of the fact that $H_*(\OZK;\k)$ is free over $\k$.
\begin{prp}[{\cite[Corollary 5.2]{gptw}}]
\label{prp:hozk_free}
If $\K$ is a flag simplicial complex, then $H_*(\OZK;\k)$ is a free $\k$-module of finite type.
\end{prp}
\begin{proof}
By the K\"unneth formula (more precisely, by the collapse of the K\"unneth spectral sequence \cite[Theorem 10.90]{rotman}), $H_*(X\times Y;\k)\simeq H_*(X;\k)\otimes H_*(Y;\k)$ if $H_*(X;\k)$ and $H_*(Y;\k)$ are free over $\k$. Hence $H_*(X\times Y;\k)$ is also a free $\k$--module.

Clearly, $H_*(S^n;\k)$ and $H_*(\Omega S^n;\k)\simeq T(a_{n-1})$ are free $\k$-modules. By Theorem \ref{thm:stanton_result} and the arguments above, the same holds for $H_*(\OZK;\k).$
\end{proof}

Hence in the flag case we have a Hopf algebra extension
$$\k\to H_*(\OZK;\k)\to H_*(\ODJ(\K);\k)\to \Lambda[m]\to\k$$
from Proposition \ref{prp:ozk_extension_statement} for any $\k.$

\section{Main calculations}
\label{sec:calculations}
In what follows, $\K$ is a flag simplicial complex on the vertex set $[m]$ with no ghost vertices, and $\k$ is a commutative ring with unit. We consider $\ZZ\times\Zm$-graded $\k$-algebras that are connected with respect to the total grading $A_n:=\bigoplus_{n=-i+|\alpha|}A_{-i,\alpha}.$
\subsection{Resolutions and formulas for differentials}
~\\
By \cite[Proposition 4.1]{cat(zk)}, the left $H_*(\ODJ(\K);\k)$-module $\k$ has a free resolution $(H_*(\ODJ(\K);\k)\otimes\kKc,d),$ $\deg\chi_\alpha:=(|\alpha|,-|\alpha|,2\alpha),$ $\deg(d)=(-1,0,0),$ with the differential
$$d(1\otimes\chi_\alpha):=\sum_{i\in\supp(\alpha)}u_i\otimes\chi_{\alpha-e_i}.$$
The isomorphism of left $H_*(\OZK;\k)$-modules
$$\Phi:H_*(\OZK;\k)\otimes\Lambda[m]\to H_*(\ODJ(\K);\k),\quad a\otimes u_I\mapsto a\cdot\widehat{u}_I$$
from Proposition \ref{prp:ozk_extension_statement} allows us to consider this resolution as a free resolution $(H_*(\OZK;\k)\otimes\Lambda[m]\otimes\kKc,\widehat{d})$ of the left  $H_*(\OZK;\k)$-module $\k.$
We apply the functor $\k\otimes_{H_*(\OZK;\k)}(-)$ and obtain a chain complex $(\Lambda[m]\otimes\kKc,\overline{d})$ whose homology is isomorphic to $\Tor^{H_*(\OZK;\k)}(\k,\k).$ The differentials $\widehat{d}$ and $\overline{d}$ are determined  by the commutative diagram
$$\xymatrix{
\dots\ar[r]
&
H_*(\ODJ(\K);\k)\otimes\kKc_{(n)}
\ar[r]^-{d}
&
H_*(\ODJ(\K);\k)\otimes\kKc_{(n-1)}
\ar[r]
&
\dots
\\
\dots\ar[r]
&
H_*(\OZK;\k)\otimes\Lambda[m]\otimes\kKc_{(n)}
\ar@{-->}[r]^-{\widehat{d}}
\ar[u]^-{\Phi\otimes\id}_-\simeq
\ar@{->>}[d]_-{\varepsilon\otimes\id\otimes\id}
&
H_*(\OZK;\k)\otimes\Lambda[m]\otimes\kKc_{(n-1)}
\ar[u]^-{\Phi\otimes\id}_-\simeq
\ar@{->>}[d]_-{\varepsilon\otimes\id\otimes\id}
\ar[r]
&
\dots
\\
\dots\ar[r]
&
\Lambda[m]\otimes\kKc_{(n)}
\ar@{-->}[r]^-{\overline{d}}
&
\Lambda[m]\otimes\kKc_{(n-1)}
\ar[r]
&
\dots
}$$
Here $\kKc_{(n)}$ is a $\k$-submodule in $\kKc$ with the basis $\{\chi_\alpha:~|\alpha|=n\}.$
With different signs, this construction was considered by the author in \cite[Section 4]{cat(zk)}. Now we describe the differential $\widehat{d}$ explicitly. For subsets $A,B\subset[m]$ define the Koszul sign $\theta(A,B):=|\{(a,b)\in A\times B:~a>b\}|.$

\begin{prp}
The differential $\widehat{d}$ is given by the formula
\begin{multline}
\label{eq:widehat_d_formula}
\widehat{d}(1\otimes u_I\otimes\chi_\alpha)=\sum_{i\in\supp(\alpha)}(-1)^{|I|}\cdot 1\otimes (u_I\wedge u_i)\otimes\chi_{\alpha-e_i}
\\
+\sum_{i\in\supp(\alpha)}
\sum_{\begin{smallmatrix}
I=A\sqcup B:\\
\max(A)>i
\end{smallmatrix}}
(-1)^{\theta(A,B)+|A|}
c(A,u_i)\otimes u_B\otimes\chi_{\alpha-e_i}.
\end{multline}
The differential $\overline{d}$ is given by the formula
\begin{equation}
\label{eq:overline_d_formula}
\overline{d}(u_I\otimes\chi_\alpha)=(-1)^{|I|}\sum_{i\in\supp(\alpha)}(u_I\wedge u_i)\otimes\chi_{\alpha-e_i}.
\end{equation}
\end{prp}
\begin{rmk}
We denote $\max(\varnothing):=-\infty$, hence $A$ cannot be empty.
\end{rmk}
\begin{proof}[Proof of the proposition.]
Recall that $u_j^2=0\in H_*(\ODJ(\K);\k).$ Therefore, by Proposition \ref{prp:u_Iu_j} we have an identity
\begin{multline*}
\widehat{u}_I\cdot u_i=
1\cdot
\begin{cases}
(-1)^{|I_{>i}|}
\widehat{u}_{I\sqcup i},
&i\notin I;\\
0,
&i\in I;
\end{cases}
+\sum_{\begin{smallmatrix}
I=A\sqcup B:\\
\max(A)>i
\end{smallmatrix}}
(-1)^{\theta(A,B)+|B|}c(A,u_i)\cdot\widehat{u}_B
\\=
\Phi\Big(1\otimes (u_I\wedge u_i)+\sum_{\begin{smallmatrix}
I=A\sqcup B:\\
\max(A)>i
\end{smallmatrix}}
(-1)^{\theta(A,B)+|B|}
c(A,u_i)\otimes u_B
\Big)\in H_*(\ODJ(\K);\k).
\end{multline*}
(Here $c(A,u_i)\in H_*(\OZK;\k)$ by Corollary \ref{crl:c(I,uj)_in_hozk}.)
Denote $\Phi_0=\Phi\otimes\id_{\kKc}.$ Then
\begin{multline*}
\Phi_0(\widehat{d}(1\otimes u_I\otimes\chi_\alpha))=d(\Phi_0(1\otimes u_I\otimes\chi_\alpha))=d(\widehat{u}_I\otimes\chi_\alpha)=(-1)^{|I|}\sum_{i\in\supp(\alpha)}\widehat{u}_Iu_i\otimes\chi_{\alpha-e_i}
\\=
(-1)^{|I|}\sum_{i\in\supp(\alpha)}
\Phi_0\Big(
1\otimes (u_I\wedge u_i)\otimes\chi_{\alpha-e_i}+
\sum_{\begin{smallmatrix}
I=A\sqcup B:\\
\max(A)>i
\end{smallmatrix}}
(-1)^{\theta(A,B)+|B|}
c(A,u_i)\otimes u_B\otimes\chi_{\alpha-e_i}
\Big).
\end{multline*}
Applying $\Phi_0^{-1}$, we obtain precisely the formula \eqref{eq:widehat_d_formula}. After the homomorphism $\varepsilon\otimes\id\otimes\id$, it turns into the formula \eqref{eq:overline_d_formula}, since $\varepsilon(1)=1$ and $\varepsilon(c(A,u_i))=0$ for $A\neq\varnothing.$
\end{proof}
\subsection{Computation of Tor-modules}
By \cite[Theorem 1.2]{cat(zk)}, for flag $\K$ we have a $\ZZ\times\ZZ\times\Zm$-graded isomorphism of $\k$-modules
\begin{equation}
\label{eq:tor_ozk_answer}
\Tor^{H_*(\OZK;\k)}(\k,\k)\cong\bigoplus_{J\subset[m]}\H_*(\K_J;\k),\quad\Tor^{H_*(\OZK;\k)}_n(\k,\k)_{-|J|,2J}\cong\H_{n-1}(\K_J;\k).
\end{equation}

Note that the homology of $\ZK$ admit a $\ZZ\times\Zm$-grading, and for any $\K$ we have a similar additive isomorphism dual to \cite[Theorem 4.5.8]{ToricTopology}:
$$H_*(\ZK;\k)\cong\bigoplus_{J\subset[m]}\H_*(\K_J;\k),\quad H_{n-|J|,2J}(\ZK;\k)\cong\H_{n-1}(\K_J;\k).$$
Hence $\Tor^{H_*(\OZK;\k)}(\k,\k)\cong H_*(\ZK;\k)$ for flag case $\K.$ Moreover, both modules are computed as the homology of $(\Lambda[u_1,\dots,u_m]\otimes\kKc,d).$

\begin{rmk}
In general, if $X$ is simply connected and $H_*(\Omega X;\k)$ is free over $\k,$ there is \emph{Milnor--Moore spectral sequence} $E^2_{p,q}=\Tor^{H_*(\Omega X;\k)}_p(\k,\k)_q\Rightarrow H_{p+q}(X;\k).$ We see that it collapses at $E^2$ for $X=\ZK$ if $\K$ is a flag complex.
For $\k=\QQ$, the collapse is explained by the \emph{coformality} of $\ZK$, see Corollary \ref{crl:zk_is_coformal} and the discussion after.
\end{rmk}

Now we construct a chain map $g$ that induces the isomorphism \eqref{eq:tor_ozk_answer}. For any chain complex $(C_\bullet,d)$ of free $\k$-modules, we have the \emph{dual complex}
$$(C^\bullet,d_{dual}),\quad C^n:=\Hom_\k(C_n,\k),\quad d_{dual}(f):c\mapsto f(d(c)).$$ Dualisation preserves isomorphisms and chain homotopies. For a simplicial complex $\K$, the augmented complex of simplicial chains $\widetilde{C}_*(\K;\k)$ has the basis $\{[I]:~I\in\K\},$ $\deg [I]:=|I|+1$ and the differential
$$d([I]):=\sum_{i\in I}(-1)^{|I_{<i}|}[I\setminus \{i\}].$$
The dual complex is the augmented complex of simplicial cochains $(\widetilde{C}^*(\K;\k),d_{dual}),$ which has the basis $\{[I]^*:~I\in\K\}$ and the differential
$$d_{dual}([I]^*)=\sum_{\begin{smallmatrix}
i\notin I:\\
I\sqcup i\in\K
\end{smallmatrix}}
(-1)^{|I_{<i}|}[I\sqcup\{i\}]^*.
$$
\begin{prp}
\label{prp:homology_hochster_map}
For every $J\subset[m],$ consider the map 
$$g_J:\widetilde{C}_{*-1}(\K_J;\k)\to (\Lambda[m]\otimes\kKc)_{*,-|J|,2J},\quad [L]\mapsto \epsilon(L,J)\cdot u_{J\setminus L}\otimes\chi_L,$$
where $\epsilon(L,J):=(-1)^{\sum_{\ell\in L}|J_{<\ell}|}.$
Then $g_J$ are chain maps, and the direct sum
$$g:\bigoplus_{J\subset[m]}\widetilde{C}_*(\K_J;\k)\to (\Lambda[m]\otimes\kKc,\overline{d})$$
induces an isomorphism on homology. Therefore,
$$H_{n,-|J|,2J}(\Lambda[m]\otimes\kKc,\overline{d})\cong\widetilde{H}_{n-1}(\K_J;\k),\quad J\subset[m],~n\geq 0,$$
all the other graded components of $H_*(\Lambda[m]\otimes\kKc,\overline{d})$ being zero.
\end{prp}
Since $\Tor^{H_*(\OZK;\k)}(\k,\k)\cong H(\Lambda[m]\otimes\kKc,\overline{d}),$ this proposition implies the formula \eqref{eq:tor_ozk_answer}. The proof is the dualisation of arguments from \cite[\S 3.2]{ToricTopology}.
\begin{proof}[Proof of Proposition \ref{prp:homology_hochster_map}]
Consider the dga algebra $(\Lambda[u_1,\dots,u_m]\otimes\kK,d)$ with the differential that is defined on generators by $d(u_i)=v_i,$ $d(v_i)=0$ and with the $\ZZ\times\ZZ\times\Zm$-grading
$$\deg u_i:=(0,-1,2e_i),\quad \deg v_i:=(1,-1,2e_i),\quad \deg d:= (1,0,0).$$

This complex has the basis $\{u_Iv^\alpha:~I\subset[m],~\alpha\in\Zm,~\supp(\alpha)\in\K\}$ and the differentials
$$d(u_I v^\alpha)=\sum_{i\in I}(-1)^{|I_{<i}|}u_{I_{<i}}v_iu_{I_{>i}}v^\alpha=\sum_{i\in I}(-1)^{|I_{<i}|}u_{I\setminus i}v_iv^\alpha.$$
Then the dual complex $(\Lambda[m]\otimes\kK)^*$ has the basis  $\{(u_Iv^\alpha)^*:~I\subset[m],~\supp(\alpha)\in\K\}$ and the differential 
$$d_{dual}((u_Iv^\alpha)^*)=\sum_{
i\in\supp(\alpha):~i\notin I
}
(-1)^{|I_{<i}|}(u_{I\sqcup i}v^{\alpha-e_i})^*.
$$ This formula is similar to \eqref{eq:overline_d_formula}. We obtain an isomorphism of chain complexes
$$\psi:(\Lambda[m]\otimes\kKc,\overline{d})\to ((\Lambda[m]\otimes\kK)^*,d_{dual}),\quad u_I\otimes\chi_\alpha\mapsto (u_Iv^\alpha)^*.$$

Consider the dga algebra $R^*(\K):=(\Lambda[m]\otimes\kK)/(u_iv_i=v_i^2=0,~i=1,\dots,m).$ It is well defined, since the ideal $(u_iv_i,v_i^2)\subset\Lambda[m]\otimes\kK$ is $d$-invariant.
The following facts are obtained in the proof of \cite[Theorem 3.2.9]{ToricTopology}.

\begin{lmm}[{\cite[Lemma 3.2.6]{ToricTopology}}]
The natural projection $\pi:\Lambda[m]\otimes\kK\to R^*(\K)$ is a chain homotopy equivalence.
\qed
\end{lmm}
\begin{lmm}
We have well defined chain maps $f_J:\widetilde{C}^*(\K_J;\k)\to R^*(\K),$
$$f_J:\widetilde{C}^{n-1}(\K_J;\k)\overset{\cong}\longrightarrow R^{n,-n,2J}(\K),\quad [L]^*\mapsto \epsilon(L,J)\cdot u_{J\setminus L}v^L,\quad \epsilon(L,J):=(-1)^{\sum_{\ell\in L}|J_{<\ell}|}.$$
The direct sum $f:\bigoplus_{J\subset[m]}\widetilde{C}^*(\K_J;\k)\to R^*(\K)$ is an isomorphism of chain complexes.
\qed
\end{lmm}
After dualisation, we obtain a chain homotopy equivalence $\pi^*$ and an isomorphism $f^*$ of chain complexes. It remains to show that the diagram
$$\xymatrix{
\bigoplus_{J\subset[m]}\widetilde{C}^*(\K_J;\k)
\ar[r]^-{g}
&
(\Lambda[m]\otimes\kKc,\overline{d})
\\
(R^*(\K))^*
\ar[r]_-{\sim}^-{\pi^*}
\ar[u]^-{f^*}_-\simeq
&
((\Lambda[m]\otimes\kK)^*,d_{dual})
\ar[u]_-{\simeq}^-\psi
}$$
is commutative. Indeed, $f^*((u_{J\setminus L}v^L)^*)=\epsilon(L,J)\cdot [L],$
hence
\[
g\Big(f^*((u_{J\setminus L}v^L)^*)\Big)=\epsilon(L,J)\cdot \epsilon(L,J) u_{J\setminus L}\otimes\chi_L = \psi\Big(\pi^*((u_{J\setminus L}v^L)^*)\Big).\qedhere
\]
\end{proof}
\begin{rmk}
In our notation, $\epsilon(L,J)=(-1)^{n},$ $n=\theta(J\setminus L,L)+|L|(|L|-1)/2$ for $L\subset J.$ 
\end{rmk}
\subsection{A chain map to the bar resolution}

\begin{thm}
\label{thm:morphism_varphi}
The identity map of the left $H_*(\OZK;\k)$-module $\k$ can be extended to the map of free resolutions
$\varphi_\bullet: (H_*(\OZK;\k)\otimes\Lambda[m]\otimes\kKc,\widehat{d})\to (\B_*(H_*(\OZK;\k)),d_{\B}),$
given by the formula
$$\varphi_n(u_I\otimes\chi_\alpha)=
(-1)^{|I|}\sum_{\begin{smallmatrix}
\alpha=e_{i_1}+\dots+e_{i_n},\\
I=A_1\sqcup\dots\sqcup A_n:\\
\max(A_t)>i_t,~\forall t\in[n]
\end{smallmatrix}}
(-1)^{\sum_{1\leq t_1<t_2\leq n}\theta(A_{t_1},A_{t_2})} \Big[c(A_1,u_{i_1})\Big|\dots\Big|c(A_n,u_{i_n})\Big].
$$
\end{thm}
\begin{proof}
We apply Corollary \ref{crl:contracting_homotopy_builds_resolution} for $\overline{\varphi}_0(u_I)=\varepsilon(u_I).$ It is sufficient to show that $\varphi_{n+1}(u_I\otimes\chi_\alpha)=s(\varphi_{n}(\widehat{d}(u_I\otimes\chi_\alpha)))$ for $|\alpha|=n+1,$ $n\geq 0.$
By \eqref{eq:widehat_d_formula} and by the $H_*(\OZK;\k)$-linearity of $\varphi_n$, we have
\begin{multline*}
\varphi_n(\widehat{d}(u_I\otimes\chi_\alpha))=\sum_{i\in\supp(\alpha)}(-1)^{|I|}\varphi_n((u_I\wedge u_i)\otimes\chi_{\alpha-e_i})
\\+
\sum_{i\in\supp(\alpha)}
\sum_{\begin{smallmatrix}
I=A\sqcup B:\\
\max(A)>i
\end{smallmatrix}}
(-1)^{\theta(A,B)+|A|}c(A,u_i)\varphi_n(u_B\otimes\chi_{\alpha-e_i}).
\end{multline*}
The map $s$ is trivial on summands of the first sum, since they belong to $\oB(H_*(\OZK;\k))\subset\Ker s\subset \B(H_*(\OZK;\k)).$ Hence we have
\begin{multline*}
s(\varphi_n(\widehat{d}(u_I\otimes\chi_\alpha)))=0+
\sum_{i\in\supp(\alpha)}
\sum_{\begin{smallmatrix}
I=A\sqcup B:\\
\max(A)>i
\end{smallmatrix}}
\sum_{\begin{smallmatrix}
\alpha-e_i=e_{i_1}+\dots+e_{i_n},\\
B=A_1\sqcup\dots\sqcup A_n:\\
\max(A_t)>i_t,~\forall t\in[n]
\end{smallmatrix}}
(-1)^{\zeta}
s\left(
c(A,u_i)\Big[c(A_1,u_{i_1})\Big|\dots\Big|c(A_n,u_{i_n})\Big]
\right)
\\=
\sum_{i\in\supp(\alpha)}
\sum_{\begin{smallmatrix}
I=A\sqcup B:\\
\max(A)>i
\end{smallmatrix}}
\sum_{\begin{smallmatrix}
\alpha-e_i=e_{i_1}+\dots+e_{i_n},\\
B=A_1\sqcup\dots\sqcup A_n:\\
\max(A_t)>i_t,~\forall t\in[n]
\end{smallmatrix}}
(-1)^{\zeta}\Big[c(A,u_i)\Big|c(A_1,u_{i_1})\Big|\dots\Big|c(A_n,u_{i_n})\Big],
\end{multline*}
where $\zeta=|B|+\theta(A,B)+|A|+\sum_{1\leq t_1<t_2\leq n}\theta(A_{t_1},A_{t_2}).$ Denoting $i=i_0,$ $A=A_0,$ we obtain
$$s(\varphi_n(\widehat{d}(u_I\otimes\chi_\alpha)))=
\sum_{\begin{smallmatrix}
\alpha=e_{i_0}+\dots+e_{i_n},\\
I=A_0\sqcup\dots\sqcup A_n:\\
\max(A_t)>i_t,~0\leq t\leq n
\end{smallmatrix}}
(-1)^{\sum_{t=1}^n\theta(A_0,A_t)+|I|+\sum_{1\leq t_1<t_2\leq n}\theta(A_{t_1},A_{t_2})}\Big[c(A_0,u_{i_0})\Big|\dots\Big|c(A_n,u_{i_n})\Big].
$$ The right hand side equals $\varphi_{n+1}(u_I\otimes\chi_\alpha)$ up to a shift of indices.
\end{proof}
\begin{thm}
\label{thm:hozk_cycles}
Let $J\subset[m].$ Let a class $\alpha\in\Tor^{H_*(\OZK)}_n(\k,\k)_{-|J|,2J}\cong \H_{n-1}(\K_J;\k)$ be represented by a cycle
$$\kappa=\sum_{I\in\K_J,|I|=n}\lambda_I\cdot [I]\in\widetilde{C}_{n-1}(\K_J;\k).$$
Then the same class is represented by the cycle $\kappa'\in\oB_n(H_*(\OZK;\k))_{-|J|,2J}$ in the bar construction,
$$\kappa':=\sum_{I\in\K_J,|I|=n}\epsilon(I,J)\lambda_I
\!\!\!\!\!\!\!
\sum_{\begin{smallmatrix}
I=\{i_1,\dots,i_n\},\\
J\setminus I = J_1\sqcup\dots\sqcup J_n:\\
\max(J_t)>i_t,~\forall t\in[n]
\end{smallmatrix}}
(-1)^{\sum_{1\leq t_1<t_2\leq n}\theta(J_{i_1},J_{i_2})}
\Big[c(J_1,u_{i_1})\Big|\dots\Big|c(J_n,u_{i_n})\Big].
$$
\end{thm}
\begin{proof}
The map $\H_*(\K_J;\k)\to\Tor^{H_*(\OZK)}_*(\k,\k)$ is induced by the composition
$$\xymatrix{
\bigoplus_{J\subset[m]}\widetilde{C}_*(\K_J;\k)\ar[r]^-g_-{\sim}
&
(\Lambda[m]\otimes\kKc,\overline{d})
\ar[r]^-{\overline{\varphi}}_\sim
&
(\oB(H_*(\OZK)),d_{\oB}),
}$$
of chain maps, where $g$ is defined in Proposition \ref{prp:homology_hochster_map} and $\overline{\varphi}$ is induced by the chain map $\varphi$ from Theorem \ref{thm:morphism_varphi}. We have $\kappa'=\overline{\varphi}(g(\kappa))$ by construction.
\end{proof}
The formulas become simpler for $n=1,2.$
\begin{crl}
\label{crl:hozk_1_cycles}
Let $J\subset[m].$ Let a class $\alpha\in\Tor^{H_*(\OZK;\k)}_2(\k,\k)_{-|J|,2J}\cong\H_1(\K_J)$ be represented by a cycle
$$\kappa=\sum_{\{i<j\}\in\K_J}\lambda_{ij}[\{i,j\}]\in\widetilde{C}_1(\K_J;\k).$$
Then the same class is represented by the following cycle in the bar construction:
\[
\kappa'=\sum_{\{i<j\}\in\K_J}(-1)^{|J_{<i}|+|J_{<j}|}
\lambda_{ij}
\!\!\!\!\!\!\!\!
\sum_{\begin{smallmatrix}
J\setminus \{i,j\}=A\sqcup B:\\
\max(A)>i,~\max(B)>j
\end{smallmatrix}}
\!\!\!\!\!\!\!\!
(-1)^{\theta(A,B)}
\Big[c(A,u_i)\Big|c(B,u_j)\Big]+(-1)^{\theta(B,A)}
\Big[c(B,u_j)\Big|c(A,u_i)\Big].
\qed
\]
\end{crl}
\begin{crl}
\label{crl_hozk_tor1_basis}
Let $J\subset[m],$ and let the simplicial complex $\K_J$ have $t+1$ path components. Let vertices $i_1,\dots,i_t,\max(J)$ be representatives of these components. Then a basis of the $\k$-module $\Tor^{H_*(\OZK;\k)}_1(\k,\k)_{-|J|,2J}\cong\H_0(\K_J;\k)\simeq\k^t$ is represented by cycles 
$$\Big[c(J\setminus i_s,u_{i_s})\Big]\in\oB_1(H_*(\OZK))_{-|J|,2J},\quad s=1,\dots,t.$$
\end{crl}
\begin{proof}
Denote $j:=\max(J).$ The cycles $\kappa_s=[\{j\}]-[\{i_s\}]\in\widetilde{C}_0(\K_J;\k),$ $1\leq s\leq t-1,$ represent a basis in $\H_0(\K_J;\k).$ By Theorem \ref{thm:hozk_cycles}, the basis in $\Tor^{H_*(\OZK;\k)}_1(\k,\k)_{-|J|,2J}$ is represented by cycles
$$
\kappa'_s=0\pm\Big[c(J\setminus i_s,u_{i_s})\Big],\quad s=1,\dots,t.
$$
(The summand $[\{j\}]$ in $\kappa_s$ does not contribute to $\kappa'_s,$ since the subset $J_1:=J\setminus\{j\}$ does not satisfy the condition $\max(J_1)>j.$)
\end{proof}
\section{Generators and relations in the flag case}
\label{sec:applications}
\subsection{Minimal sets of generators}
Denote $\b_0(X):=\rank\H_0(X;\k).$ This number does not depend on $\k,$ since $\b_0(X)+1$ is the number of path components in $X$.

\begin{thm}
\label{thm:minimal_generators}
Let $\K$ be a flag simplicial complex on vertex set $[m]$ and $\k$ be a commutative ring with unit. For every $J\subset[m],$ choose a $\b_0(\K_J)$-element subset $\Theta(J)\subset J\setminus\{\max(J)\}$ such that $\Theta(J)\sqcup\{\max(J)\}$ contains exactly one vertex from each path component of $\K_J.$ Then $H_*(\OZK;\k)$ is multiplicatively generated by the following set of $\sum_{J\subset[m]}\b_0(\K_J)$ elements:
$$\Big\{c(J\setminus i,u_i):\quad i\in\Theta(J),~J\subset[m]\Big\},\quad c(J\setminus i,u_i)\in H_{-|J|,2J}(\OZK;\k).$$
If $\k$ is a principal ideal domain, this set is minimal: any $\ZZ\times\Zm$-homogeneous presentation of $H_*(\OZK;\k)$ contains at least $\b_0(\K_J)$ generators of degree $(-|J|,2J);$ any $\ZZ$-homogeneous presentation contains at least $\sum_{|J|=n}\b_0(\K_J)$ generators of degree $n.$
\end{thm}
\begin{proof}
By Corollary~\ref{crl_hozk_tor1_basis}, images of cycles $\{[c(J\setminus i,u_i)]:J\subset[m],i\in\Theta(J)\}\subset\oB_1(H_*(\OZK;\k))$ additively generate the $\k$-module $\Tor^{H_*(\OZK;\k)}_1(\k,\k).$ Hence, by Theorem~\ref{thm:presentation_from_cycles}(1), the algebra $H_*(\OZK;\k)$ is multiplicatively generated by the elements in question.
The lower bounds on the number of generators follow from the formula \eqref{eq:tor_ozk_answer} and from Theorem~\ref{thm:size_of_pres_PID}.(2).
\end{proof}
\begin{dfn}
\label{dfn:GPTW_gens}
Let $\K$ be a simplicial complex on $[m],$ and let $J\subset[m].$ Choose $\Theta(J)$ as the set of the smallest vertices in corresponding path components. More precisely, define $\Theta(J)$ as the set of all vertices $i\in J$ such that
\begin{enumerate}
\item $i$ and $\max(J)$ belong to different path components of the complex $\K_J;$
\item $i$ is the smallest vertex (has the smallest number) in its path component.
\end{enumerate}
The corresponding set of generators $\{c(J\setminus i,u_i):~i\in\Theta(J),~J\subset[m]\}$ will be called the \emph{GPTW generators}.
\end{dfn}
Grbi\'c, Panov, Theriault and Wu proved \cite[Theorem 4.3]{gptw} that GPTW generators minimally generate the algebra $H_*(\OZK;\k)$ if $\k$ is a field. The minimality was proved using topological methods. Our Theorem~\ref{thm:minimal_generators} gives a purely algebraic proof for any ring $\k$.
\subsection{Rewriting of nested commutators}
Thus the GPTW generators are indeed multiplicative generators of the algebra $H_*(\OZK;\k)$ for any ring $\k$ and any flag complex $\K$.
\begin{dfn}
Let $i\in J\subset[m].$ Express the element $c(J\setminus i,u_i)\in H_*(\OZK;\k)$ as a non-commutative polynomial in GPTW generators (this expression may be non-unique). Any such expression will be denoted by $\widehat{c}(J\setminus i,u_i).$
\end{dfn}

These non-commutative polynomials can be computed recursively, following the proof of \cite[Theorem 4.3]{gptw}. We describe an explicit rewriting process. 

\begin{algo}
\label{algo:rewriting}
Suppose that expressions $\widehat{c}(A\setminus t,u_t),$ $|A|<|J|,$ are already computed, and we compute $\widehat{c}(J\setminus i,u_i)$.
Three cases is possible:
\begin{enumerate}
\item $i=\max(J).$ Denote $j=\max(J\setminus i).$ Then $c(J\setminus i,u_i)=c(J\setminus ij,[u_j,u_i])=c(J\setminus j,u_j).$ The task is reduced to the case $i\neq\max(J).$ 
\item $i$ and $\max(J)$ belong to the same path component of $\K_J$. The length of the shortest path from $i$ to $\max(J)$  along the edges of $\K_J$ will be called the \emph{rank} of a vertex $i.$ We proceed by induction on the rank. The case of rank zero is discussed above. If rank equals 1, we have $[u_{\max(J)},u_i]=0,$ so $$c(J\setminus i,u_i)=c(J\setminus\{i,\max(J)\},[u_{\max(J)},u_i])=0.$$ Suppose that rank is greater than one, and let $\{i,j\}$ be the first edge in (any) shortest path from $i$ to $\max(J).$ Since $[u_i,u_j]=0\in H_*(\OZK;\k),$ the identity \eqref{eq:commutator_hard_identity} expresses $c(J\setminus i,u_i)$ in terms of $c(J\setminus j,u_j)$ (this element has smaller rank) and commutators of smaller degree (expressions for which are already computed).
\item $i$ and $\max(J)$ are in different path components. Let $i_0$ be the smallest vertex of the component that contains $i.$ The length of the shortest path from $i$ to $i_0$ will be called the \emph{rank} of a vertex $i.$ If the rank is zero, then $i\in\Theta(J),$ so we can set $\widehat{c}(J\setminus i,u_i):=c(J\setminus i,u_i).$ Otherwise we decrease the rank using \eqref{eq:commutator_hard_identity}, as in case (2).
\end{enumerate}
\end{algo}
\begin{rmk}
Similar argument works more generally: suppose that we have a set of elements $\{x_{J,i}:i\in J\subset [m]\}$ such that, for any $\{i,j\}\in\K_J,$ the linear combination $x_{J,i}\pm x_{J,j}$ is a non-commutative polynomial on elements of smaller degree. Then we can express each element $x_{A,t}$ throught the ``GPTW elements'' $\{x_{J,i}:i\in\Theta(J),J\subset[m]\}$ by a similar rewriting process. In our case $x_{J,i}=c(J\setminus i,u_i),$ and the polynomial is given by the last summand in \eqref{eq:commutator_hard_identity}.
\end{rmk}
\subsection{Minimal sets of relations}
Let $M$ be a finitely generated $\k$-module. Denote the smallest number of generators by $\gen(M).$ Denote $b_0(X):=\gen(\H_0(X;\k)),$ $b_1(X;\k):=\gen(H_1(X;\k)).$
\begin{thm}
\label{thm:minimal_relations}
Let $\K$ be a flag simplicial complex on vertex set $[m]$, $\k$ be a commutative ring. For each $J\subset[m]$, choose a collection of simplicial $1$-cycles
$$\sum_{\{i<j\}\in\K_J}\lambda_{ij}^{(\alpha)}[\{i,j\}]\in\widetilde{C}_1(\K_J;\k)$$
that generate the $\k$-module $H_1(\K_J;\k).$ Then the algebra $H_*(\OZK;\k)$ is presented by GPTW generators $\{c(J\setminus i,u_i):~i\in\Theta(J),~J\subset[m]\}$ (see Definition \ref{dfn:GPTW_gens}) modulo the relations
\begin{equation}
\label{eq:relations_hozk}
\sum_{\{i<j\}\in\K_J}(-1)^{|J_{<i}|+|J_{<j}|}\lambda_{ij}^{(\alpha)}\sum_{\begin{smallmatrix}
J\setminus\{i,j\}=A\sqcup B:\\
\max(A)>i,~\max(B)>j
\end{smallmatrix}}
(-1)^{\theta(A,B)+|A|}
\Big[\widehat{c}(A,u_i),\widehat{c}(B,u_j)\Big]=0.
\end{equation}
In particular, $H_*(\OZK;\k)$ admits a presentation by $\sum_{J\subset[m]}\b_0(\K_J)$ generators modulo $\sum_{J\subset[m]}b_1(\K_J;\k)$ relations: one should take the 1-cycles that correspond to minimal sets of generators.

If $\k$ is a principal ideal domain, this presentation is minimal: any  $\ZZ\times\Zm$-homogeneous presentation contain at least $b_1(\K_J;\k)$ relations of degree $(-|J|,2J)$ for every $J\subset[m].$
\end{thm}
\begin{proof}
By Corollary~\ref{crl:hozk_1_cycles}, our 1-cycles correspond to the elements
$$\sum_{\{i<j\}\in\K_J}(-1)^{|J_{<i}|+|J_{<j}|}\lambda_{ij}^{(\alpha)}\sum_{\begin{smallmatrix}
J\setminus ij = A\sqcup B:\\
\max(A)>i,~\max(B)>j
\end{smallmatrix}}
(-1)^{\theta(A,B)} \Big[c(A,u_i)\Big|c(B,u_j)\Big]+(-1)^{\theta(B,A)}\Big[c(B,u_j)\Big|c(A,u_i)\Big]$$
in bar construction, and their images additively generate $\Tor^{H_*(\OZK;\k)}_2(\k,\k).$
We apply Theorem~\ref{thm:presentation_from_cycles}(2) to this situation. (In the notation of this theorem, we take GPTW generators as $a_1,\dots,a_N.$ Their images freely generate $\Tor^{H_*(\OZK)}_1(\k,\k),$ so we can take $R=0.$ We take $\widehat{c}(A,u_i)$ and $\widehat{c}(B,u_j)$ as polynomials $P_{j,\alpha}$ and $Q_{j,\alpha}$.) It follows that $H_*(\OZK;\k)$ is generated by GPTW generators and presented by the relations
$$\sum_{\{i<j\}\in\K_J}(-1)^{|J_{<i}|+|J_{<j}|}\lambda_{ij}^{(\alpha)}\sum_{\begin{smallmatrix}
J\setminus ij = A\sqcup B:\\
\max(A)>i,~\max(B)>j
\end{smallmatrix}}
(-1)^{\theta(A,B)} \overline{\widehat{c}(A,u_i)}\widehat{c}(B,u_j)+(-1)^{\theta(B,A)}\overline{\widehat{c}(B,u_j)}\widehat{c}(A,u_i) = 0.$$
Denote $x=\widehat{c}(A,u_i),$ $y=\widehat{c}(B,u_j).$ Since $\theta(A,B)+\theta(B,A)\equiv |A|\cdot |B|,$ we have
\begin{multline*}
(-1)^{\theta(A,B)}\overline{x}y+(-1)^{\theta(B,A)}\overline{y}x=(-1)^{\theta(A,B)}\left((-1)^{|A|}xy+(-1)^{|B|+|A|\cdot|B|}yx\right)
\\=
(-1)^{\theta(A,B)+|A|}\left(
xy-(-1)^{(|A|+1)(|B|+1)}yx
\right)=(-1)^{\theta(A,B)+|A|}[x,y].
\end{multline*}
Hence the obtained relations coincide with \eqref{eq:relations_hozk}. Finally, the lower bound on the number of relations follows from \eqref{eq:tor_ozk_answer} and Theorem~\ref{thm:size_of_pres_PID}(2).
\end{proof}
Sometimes we can reduce the number of relations if the presentation is not required to be $\ZZ\times\Zm$-homogeneous. For example, suppose that for some $I,J\subset[m]$ we have $|I|=|J|=n,$ $H_1(\K_I;\ZZ)=\ZZ/2,$ $H_1(\K_J;\ZZ)=\ZZ/3.$ Then the graded components of the module $\Tor^{H_*(\OZK;\ZZ)}_2(\ZZ,\ZZ)$ having multidegrees $(-n,2I)$ and $(-n,2J)$ are equal to $\ZZ/2$ and $\ZZ/3.$
By Theorem~\ref{thm:size_of_pres_PID}, every $\ZZ\times\Zm$-homogeneous presentation of $H_*(\OZK;\ZZ)$ should contain relations of these multidegrees. On the other hand, these $\ZZ\times\Zm$-graded components contribute $\ZZ/2\oplus\ZZ/3\simeq\ZZ/6$ to the $\ZZ$-graded component of degree $n.$ Hence we can take just one $\ZZ$-homogeneous relation (for example, the sum of these $\ZZ\times\Zm$-homogeneous relations). Let us give a general result. 
\begin{thm}
\label{thm:hozk_z_homogeneous}
Let $\K$ be a flag simplicial complex and $\k$ be a principal ideal domain. Consider all homogeneous presentations of the $\ZZ$-graded $\k$-algebra $H_*(\OZK;\k).$
\begin{enumerate}
\item There is a presentation that consists of, for each $n\geq 0$, exactly $\sum_{|J|=n}\b_0(\K_J)$ generators and exactly $\gen(\bigoplus_{|J|=n}H_1(\K_J;\k))$ relations of degree $n.$ One can take GPTW generators as generators,
and take linear combinations of identities from Theorem~\ref{thm:minimal_relations}, corresponding to minimal generators of the $\k$-module $\bigoplus_{|J|=n}H_1(\K_J;\k),$ as relations.
\item For every $n\geq 0$, any presentation contains at least $\sum_{|J|=n}\b_0(\K_J)$ generators and at least $\gen(\bigoplus_{|J|=n}H_1(\K_J;\k))$ relations of degree $n.$
\end{enumerate}
\end{thm}
\begin{proof}
By Theorem~\ref{thm:size_of_pres_PID}, the number $\gen\Tor^{H_*(\OZK;\k)}_1(\k,\k)_n$ (the number $\gen\Tor^{H_*(\OZK;\k)}_2(\k,\k)_n+\rel\Tor^{H_*(\OZK;\k)}_1(\k,\k)_n$) is a precise bound on the number of generators (of relations) of degree $n.$ By~\eqref{eq:tor_ozk_answer}, we have
$$\Tor^{H_*(\OZK;\k)}_1(\k,\k)_n=\bigoplus_{|J|=n}\H_0(\K_J;\k)\simeq\k^{\sum_{|J|=n}\b_0(\K_J)},~\Tor^{H_*(\OZK;\k)}_2(\k,\k)=\bigoplus_{|J|=n}H_1(\K_J;\k);$$
hence $\gen\Tor_1=\sum_{|J|=n}\b_0(\K_J)$ and $\rel\Tor_1=0.$
One can take the GPTW generators since the images of corresponding cycles generate $\Tor_1$ by Corollary~\ref{crl_hozk_tor1_basis}.
\end{proof}

\subsection{Example: moment-angle complexes for $m$-cycles}
Let $\K$ be the boundary of $m$-gon. The corresponding moment-angle complex $\ZK$ is homeomorphic to a connected sum of sphere products, $\ZK\cong \#_{k=3}^{m-1} (S^k\times S^{m+2-k})^{\#(k-2)\binom{m-2}{k-1}},$
and hence $H_*(\OZK;\k)$ is a one-relator algebra. It was considered in \cite{verevkin,onerelator}. From the point of view of Theorem~\ref{thm:minimal_relations}, the relation corresponds to the 1-cycle
$$\kappa=[\{1,m\}]-\sum_{i=1}^{m-1}[\{i,i+1\}]\in\widetilde{C}_1(\K;\k)$$
and has the form
$$
\sum_{\begin{smallmatrix}
\{2,\dots,m-1\}=A\sqcup B:\\
\max(A)>1,~\max(B)>m
\end{smallmatrix}}
\!\!\!\!\!\!\!\!
(\dots)
\quad
-
\quad
\sum_{i=1}^{m-1}(-1)^{(i-1)+i}
\!\!\!\!\!\!\!\!
\sum_{\begin{smallmatrix}
[m]\setminus\{i,i+1\}=A\sqcup B:\\
\max(A)>i,~\max(B)>i+1
\end{smallmatrix}}
(-1)^{\theta(A,B)+|A|}\Big[\widehat{c}(A,u_i),\widehat{c}(B,u_{i+1})\Big]=0.$$
The first sum is empty, since $\max(B)\leq m-1.$ Similarly, in the second sum the inner sum is empty for $i=m-1,m-2.$ The simplified relation is
$$\sum_{i=1}^{m-3}\sum_{
\begin{smallmatrix}
[m]\setminus\{i,i+1\}=A\sqcup B:\\
\max(A),\max(B)\geq i+2
\end{smallmatrix}}
(-1)^{\theta(A,B)+|A|}\Big[\widehat{c}(A,u_i),\widehat{c}(B,u_{i+1})\Big]=0.$$
Some summands are immediately zero. For example, if $\max(B)=i+2,$ then $c(B,u_{i+1})=c(B\setminus\{i+2\},[u_{i+2},u_{i+1}])=0,$ so we can take $\widehat{c}(B,u_{i+1})=0.$ Similarly, $c(A,u_1)=0$ if $i=1$ and $m\in A.$ Other summands can be computed using Algorithm~\ref{algo:rewriting}. We were not able to obtain a closed formula for this relation (as a polynomial of GPTW generators or other minimal generators).
However, we at least have an effective algorithm that computes the relation for any given $m.$

Consider the case $m=5.$ Besides from the partitions $[5]\setminus\{i,i+1\}=A\sqcup B$ considered above, for $i=1$ the allowed partitions are
$\{3,4,5\}=\{3\}\sqcup\{4,5\}=\{4\}\sqcup\{3,5\}=\{3,4\}\sqcup\{5\};$
for $i=2$ the allowed partitions are $\{1,4,5\}=\{4\}\sqcup\{1,5\}=\{1,4\}\sqcup\{5\}.$
The resulting relation has five summands:
\begin{multline*}
(-1)^{\theta(3,45)+1}
\Big[\widehat{c}(3,u_1),\widehat{c}(45,2)\Big]+
(-1)^{\theta(4,35)+1}
\Big[\widehat{c}(4,u_1),\widehat{c}(35,2)\Big]+
(-1)^{\theta(34,5)+2}
\Big[\widehat{c}(34,u_1),\widehat{c}(5,u_2)\Big]
\\+
(-1)^{\theta(4,15)+1}
\Big[\widehat{c}(4,u_2),\widehat{c}(15,u_3)\Big]+
(-1)^{\theta(14,5)+2}
\Big[\widehat{c}(14,u_2),\widehat{c}(5,u_3)\Big]=0.
\end{multline*}
All commutators, apart from $\widehat{c}(14,u_2)=[u_1,[u_4,u_2]]=-[u_2,[u_4,u_1]]=-c(24,u_1),$ already are GPTW generators. We obtain the following identity between the generators:
\begin{multline*}
-\Big[[u_3,u_1],[u_4,[u_5,u_2]]\Big]+\Big[[u_4,u_1],[u_3,[u_5,u_2]]\Big]-\Big[[u_5,u_2],[u_3,[u_4,u_1]]\Big]
\\+
\Big[[u_4,u_2],[u_1,[u_5,u_3]]\Big]+\Big[[u_5,u_3],[u_2,[u_4,u_1]]\Big]=0.
\end{multline*}
This relation was first obtained by Veryovkin as a result of bruteforce \cite[Theorem 3.2]{verevkin}. For $m=6$, the analogous relation is initially the sum of $7+10+4=21$ commutators. After computing the elements $\widehat{c}(J\setminus i,u_i)$ and changing the set of generators, it can be written as $\sum_{i=1}^{17}[a_i,b_i]=0$ (see \cite[Theorem 4.1]{verevkin}). This agrees with the homeomorphism $\ZK\cong (S^3\times S^5)^{\# 9}\# (S^4\times S^4)^{\# 8}.$

\section{Homotopical properties in the flag case}
\label{sec:homotopy}
\subsection{Homotopy groups}
As in \cite{stanton}, we denote by $\mathcal{P}$ the class of H-spaces which are homotopy equivalent to finite type products of spheres and loops on simply connected spheres, and by $\mathcal{W}$ the class of topological spaces which are homotopy equivalent to finite type wedges of simply connected spheres. The author thanks Lewis Stanton for providing a proof of the following lemma.

\begin{lmm}
\label{lmm:eliminating_spheres}
Let $A_1,\dots,A_m$ be connected topological spaces, $\K$ be a simplicial complex on $[m]$, and suppose that $\Omega(\underline{CA},\underline{A})^\K\in\mathcal{P}$. Then $\Omega(\underline{CA},\underline{A})^\K$ is homotopy equivalent to a finite type product of loops on simply connected spheres.
\end{lmm}
\begin{proof}
By \cite[Corollary 9.8]{theriault_toric}, $\Omega(\underline{CA},\underline{A})^\K\simeq\prod_{i=1}^m\Omega\Sigma Y_i$ for some spaces $Y_i.$ Since the class $\mathcal{P}$ is closed under retracts \cite[Theorem 3.10]{stanton}, $\Omega\Sigma Y_i\in\mathcal{P}.$ By repeated use of the homotopy equivalence $\Sigma(X\times Y)\simeq\Sigma X\vee\Sigma Y\vee \Sigma(X\wedge Y)$ and the James splitting $\Sigma\Omega\Sigma X\simeq\bigvee_{n\geq 1}\Sigma X^{\wedge n}$ \cite{james}, we have $\Sigma Z\in\mathcal{W}$ for $Z\in\mathcal{P}$. In particular, $\Sigma\Omega\Sigma Y_i\in\mathcal{W}$. On the other hand, $\Sigma Y_i$ is a retract of $\Sigma\Omega\Sigma Y_i$ by the James splitting. The class $\mathcal{W}$ is closed under retracts (see for example \cite[Lemma 3.1]{amelotte}), so $\Sigma Y_i\in\mathcal{W}.$ Now $\Omega\Sigma Y_i$ is homotopy equivalent to a product of loops on spheres by the Hilton--Milnor theorem. It follows that the same holds for $\prod_{i=1}^m\Omega \Sigma Y_i.$
\end{proof}

\begin{proof}[Proof of Theorem~\ref{thm:homotopy_groups_flag}]
Since $\K$ is flag, we have $\Omega\ZK\in\mathcal{P}$ by Theorem~\ref{thm:stanton_result}. Hence $\Omega\ZK=\Omega(C S^1,S^1)^\K$ is a product of loops on spheres by Lemma \ref{lmm:eliminating_spheres}. It follows that for some $D_n\geq 0$ we have a homotopy equivalence
$$\OZK\simeq \prod_{n\geq 2}(\Omega S^n)^{\times D_n}.$$
The numbers $D_n$ are finite, since $\dim_\k H_i(\OZK;\k)<\infty$ for all $i.$ (Here $\k$ is any field.)
Also $D_2=0,$ since $\ZK$ is 2-connected \cite[Proposition 4.3.5]{ToricTopology}. In order to compute $D_n,$ we calculate $\dim H_i(\OZK;\k)$ twice. Recall that the \emph{Poincar\'e series} $P(V;t)$ of a graded $\k$-vector space $V$ are the formal power series
$$P(V;t):=\sum_{i\geq 0}\dim_\k(V_i)\cdot t^i\in\ZZ[[t]].$$ We have $P(V\oplus W)=P(V;t)+P(W;t)$ and $P(V\otimes W;t)=P(V;t)\cdot P(W;t).$
From $F(H_*(\Omega S^k;\k);t)=(1-t^{k-1})^{-1}$ and the K\"unneth formula we have $F(H_*(\OZK;\k);t)=\prod_{n\geq 3}(1-t^{n-1})^{-D_n}.$
On the other hand, it is known (see \cite[Proposition 8.5.4]{ToricTopology} and \cite[Theorem 4.8]{cat(zk)}) that
$$F(H_*(\OZK;\k);t)=\frac{1}{(1+t)^{m-d}\cdot h_\K(-t)}=-\frac{1}{\sum_{J\subset[m]}\widetilde{\chi}(\K_J)t^{|J|}}$$
for a flag complex $\K.$ We obtain the required identity \eqref{eq:BCD_n}.
\end{proof}
\begin{rmk}
In the proof above, the algebra $H_*(\OZK;\k)$ is actually $\ZZ\times\Zm$-graded. We expect that factors of the product \eqref{eq:OZK_homotopy_type} can be considered as ``$\ZZ\times\Zm$-graded spheres'', and thus $\pi_*(\OZK)$ admits a functorial $\ZZ\times\Zm$-grading as conjectured in \cite[Remark 4.10]{cat(zk)}.
\end{rmk}
\begin{prb}
Describe the Whitehead bracket in $\pi_*(\ZK)$ in terms of the decomposition \eqref{eq:ZK_homotopy_ranks}.
\end{prb}

\subsection{Rational coformality of moment-angle complexes}
Let $X$ be a simply connected space and $\Omega X$ be the space of Moore loops. Since $\Omega X$ is a strictly associative topological monoid, the chain complex $C_*(\Omega X;\k)$ is a dga algebra with respect to the Pontryagin product for any $\k$. Also, the cochain complex $C^*(X;\k)$ is a dga algebra with respect to the Kolmogorov-Alexander product (cup product).
\begin{dfn}
A topological space $X$ is \emph{formal} over a ring $\k$, if the dga algebras $H^*(X;\k)$ (with zero differential) and $C^*(X;\k)$ are quasi-isomorphic (are connected by a zigzag of dga maps which induce isomorphisms on homology).
\end{dfn}
\begin{dfn}
A simply connected space $X$ is \emph{coformal} over a ring $\k$, if the dga algebras $H_*(\Omega X;\k)$ (with zero differential) and $C_*(\Omega X;\k)$ are quasi-isomorphic.
\end{dfn}
The notions of formality and coformality (over a field of characteristic zero) arose in rational homotopy theory, and were initially formulated in terms of Sullivan and Quillen models. The rational homotopy type of a formal (coformal) space is fully determined by the algebra $H^*(X;\QQ)$ (by the algebra $H_*(\Omega X;\QQ)$). As proved by Saleh \cite[Corollary 1.2, 1.4]{saleh}, our definitions are equivalent to the classical ones.

It is known \cite[Theorem 4.8]{notbohm_ray} that all Davis--Januszkiewicz spaces $\DJ(\K)$ are formal over $\ZZ$ (therefore, over any ring $\k$). Also, $\DJ(\K)$ is coformal over $\QQ$ if and only if $\K$ is flag \cite[Theorem 8.5.6]{ToricTopology}. First examples of non-formal moment-angle complexes were constructed by Baskakov \cite{baskakov} using Massey products. See \cite[Introduction]{buchstaber_limonchenko} for a survey of further developments in this area.

The following result of Huang can be used to prove coformality over $\QQ.$
\begin{prp}[{\cite[Proposition 5.1]{huang}}]
\label{prp:huang}
Let $F\overset{i}\longrightarrow E\to B$ be a fibration of nilpotent spaces of finite type, such that
\begin{itemize}
\item The map $i_*:\pi_*(F)\otimes_\ZZ\QQ\to\pi_*(E)\otimes_\ZZ\QQ$ is injective;
\item $E$ is coformal over $\QQ.$
\end{itemize}
Then $F$ is coformal over $\QQ.$\qed
\end{prp}
\begin{crl}
\label{crl:zk_is_coformal}
Let $\K$ be a flag simplicial complex with no ghost vertices. Then $\ZK$ is coformal over $\QQ.$
\end{crl}
\begin{proof}
We apply Proposition \ref{prp:huang} to the fibration $\ZK\to\DJ(\K)\to B\TT^m.$ By Proposition~\ref{prp:zk_dj_bt_fibration} and exact sequence of homotopy groups, $\pi_*(\ZK)\to\pi_*(\DJ(\K))$ is injective. The second condition holds by \cite[Theorem 8.5.6]{ToricTopology}.
\end{proof}
It is natural to hope that Huang's theorem admits the following generalisation.
\begin{cnj}
Let $F\to E\overset{p}\longrightarrow B$ be a fibration of simply connected spaces of finite type, such that
\begin{itemize}
\item $\Omega p$ has a homotopy section;
\item $E$ is coformal over $\k.$
\end{itemize}
Then $F$ is coformal over $\k.$
\end{cnj}

Let $X$ be a simply connected space such that $H_*(\Omega X;\k)$ is a free $\k$-module. The tensor filtration on the bar construction $\oB(C_*(\Omega X;\k))$ gives rise to the \emph{Milnor--Moore spectral sequence}
$$E^2_{p,q}=\Tor^{H_*(\Omega X;\k)}_p(\k,\k)_q\Rightarrow Tor^{C_*(\Omega X;\k)}_{p+q}(\k,\k)\cong H_{p+q}(X;\k).$$
(The last isomorphism is due to Eilenberg--Moore, see \cite[Theorem IV]{adams_cobar}).

The differential $Tor$ is preserved by quasi-isomorphisms. Hence the spectral sequence collapses at $E^2$ if $X$ is coformal over $\k$. On the other hand, it collapses for $X=\ZK$ in the flag case, see \eqref{eq:tor_ozk_answer}. This suggests the following conjecture.
\begin{cnj}
Let $\K$ be a flag simplicial complex. Then the spaces $\DJ(\K)$ and $\ZK$ are coformal over any commutative ring with unit.
\end{cnj}

\subsection{A necessary condition for the rational formality in the flag case}
The space $X$ is \emph{Koszul} if it is both formal and coformal over $\QQ.$ Hence $\DJ(\K)$ is Koszul if and only if $\K$ is flag. Koszul spaces were introduced by Berglund \cite{berglund_koszul}.

\begin{dfn}
Let $\k$ be a field, $A=\bigoplus_{n\in\ZZ}A^n$  be a graded $\k$-algebra that admits an additional ``weight'' grading
$A^n=\bigoplus_{j\geq 0}A^{n,(j)}.$ The algebra $A$ is \emph{Koszul} with respect to the weight grading if $\Ext_A^i(\k,\k)^{n,(j)}=0$ for all $i\neq j.$
\end{dfn}
For every Koszul algebra, there is a \emph{quadratic dual} Koszul algebra $A^!,$ see \cite{froeberg_koszul}. More explicitly, we set
$$A^!:=\Ext_A(\k,\k),\quad (A^!)^{n,(i)}=\Ext^i_A(\k,\k)^{-i-n,(i)}.$$
Then it is known that $(A^!)^!\cong A$ as bigraded algebras.

\begin{rmk}
In the classical theory of Koszul algebras \cite{priddy,froeberg_koszul} the $\ZZ$-grading $A=\bigoplus_{n\in\ZZ}A^n$ is absent, and only the weight grading $(A^!)^{(i)}=\Ext^i_A(\k,\k)^{(i)}$ is considered. Classical results are readily generalised to the graded case.
\end{rmk}

The following result is due to Berglund. Note that we replace the Koszul Lie algebra with their universal enveloping algebras. Berglund considers a stonger version of the Koszul duality, the duality between Lie algebras and commutative algebras.
\begin{thm}[{\cite[Theorem 2, Theorem 3]{berglund_koszul}}]
\label{thm:from_berglund}
Let $X$ be a simply connected space of finite type such that $X$ is coformal over $\QQ.$ The following conditions are equivalent:
\begin{enumerate}
\item[(a)] $X$ is formal over $\QQ;$
\item[(b)] The graded algebra $A=H_*(\Omega X;\QQ)$ admits a weight grading $A=\bigoplus_{i\geq 0}A^{(i)}$ such that $A$ is Koszul with respect to it.
\end{enumerate}
Moreover, if these conditions are met, then the $\ZZ$-graded algebras $A^!$ and $H^{-*}(X;\QQ)$ are isomorphic: $H^n(X;\QQ)\cong\bigoplus_{i\geq 0}(A^!)^{-n,(i)}.$\qed
\end{thm}
\begin{thm}
\label{thm:zk_flag_formality_criterion}
Let $\K$ be a flag simplicial complex on $[m]$ with no ghost vertices, such that $\ZK$ is rationally formal. Then $\Gamma=H^*(\ZK;\QQ)$ is a Koszul algebra with respect to the grading
$$\Gamma^{(i)}:=\bigoplus_{J\subset[m]} H^{i-|J|,2J}(\ZK;\QQ)=\bigoplus_{J\subset[m]}\H^{i-1}(\K_J;\QQ).$$
In particular, $\Gamma$ is generated by elements in $\H^0(\K_J;\QQ)$ modulo the relations in $\H^1(\K_J;\QQ).$
\end{thm}
\begin{proof}
By Theorem~\ref{thm:from_berglund}, the algebra $A=H_*(\OZK;\QQ)$ is Koszul with respect to a weight grading $A=\bigoplus_{i\geq 0}A^{(i)}.$ From \cite[Theorem 1.2]{cat(zk)} we have $\Tor^A_i(\QQ,\QQ)_j=\bigoplus_{|J|=j}\H_{i-1}(\K_J;\QQ)$. Therefore, $\Ext_{A}^i(\QQ,\QQ)^{j}=\bigoplus_{|J|=j}\H^{i-1}(\K_J;\QQ).$
The algebra $A$ is Koszul, hence $$\Ext_A^{i}(\QQ,\QQ)^j=\Ext_A^{i}(\QQ,\QQ)^{j,(i)}=(A^!)^{-i-j,(i)}.$$ Since $(A^!)^*\cong \Gamma^{-*}$ as graded algebras, we obtain a weight grading
$$\Gamma^{i+j,(i)}=\bigoplus_{|J|=j}\H^{i-1}(\K_J;\QQ),\quad \Gamma^{(i)}=\bigoplus_{J\subset[m]}\H^{i-1}(\K_J;\QQ)$$
such that $\Gamma$ is Koszul with respect to it. Finally, any Koszul algebra is generated by elements of weight 1 modulo relations of weight 2. 
\end{proof}
\begin{cnj}
If $\K$ is flag and $H^*(\ZK;\QQ)$ is Koszul with respect to the grading from Theorem~\ref{thm:zk_flag_formality_criterion}, then $\ZK$ is formal over $\QQ.$
\end{cnj}
\appendix
\section{Presentations of connected graded algebras}
\label{sec:appendix_presentations}
In this section we prove Theorems \ref{prp:presentation_exact_sequence} and \ref{thm:size_of_pres_PID} that generalise some results of Wall \cite[Section 7]{wall}. We also prove Theorem \ref{thm:presentation_from_cycles}, which seems to be new. We use the notations from Section \ref{sec:algebra}; some of them are recalled below.
\subsection{Conventions}
The ring $\k$ is assumed to be an arbitrary commutative associative ring with unit. All tensor products are over $\k.$

We consider $G$-graded $\k$-algebras, where $G$ is a commutative monoid supplied with a homomorphism $G\to\ZZ.$ It induces a $\ZZ$-grading. Such algebra $A$ is \emph{connected} if it is connected with respect to the $\ZZ$-grading, i.e. $A_{<0}=0$ and $A_0=\k\cdot 1.$ Then the standard augmentation $\varepsilon:A\to A_0\cong\k$ makes $\k$ a left $A$-module and a right $A$-module.

Every complex of $G$-graded modules is considered as a $\ZZ\times G$-graded module with a differential of degree $(-1,0).$ Hence, $A$-linear differentials satisfy the following version of Leibniz's rule:
$$d(a\cdot x)=(-1)^{\deg(a)}a\cdot d(x)=-\overline{a}\cdot d(x),$$
where $\overline{a}:=(-1)^{1+\deg(a)}a.$

A \emph{presentation} of a connected $\k$-algebra $A$ is an isomorphism of the form $A\simeq T(x_1,\dots,x_N)/(r_1,\dots,r_M),$ sometimes written as
$$A\simeq T(x_1,\dots,x_N)/(r_1=\dots=r_M=0),$$
where $T(x_1,\dots,x_N)$ is a tensor algebra and $(r_1,\dots,r_M)\subset T(x_1,\dots,x_N)$ is the two-sided ideal generated by the set $\{r_1,\dots,r_M\}.$ It is assumed that generators and relations are homogeneous and have positive degree, hence belong to $\Ker\varepsilon.$ Note that $A$ is not required to be a free $\k$-module, and $M,$ $N$ can be infinite of any cardinality.

\subsection{Exact sequence of a presentation}
Let $T(x_1,\dots,x_N)$ be a tensor algebra generated by homogeneous elements of positive degrees. Every element $w\in T(x_1,\dots,x_N)$ is uniquely represented as a sum
$$w=\varepsilon(w)+\sum_{i=1}^N w_i\cdot x_i,\quad w_i\in T(x_1,\dots,x_N).$$
In the next proposition we use this representation implicitly. For example, we assume that $r_j =\varepsilon(r_j)+\sum_{i=1}^N r_{ji}\cdot x_i.$ Since $r_j\in\Ker\varepsilon,$ the first summand is zero.
\begin{prp}
\label{prp:presentation_exact_sequence}
Let $A=T(x_1,\dots,x_N)/(r_1,\dots,r_M)$ be a presentation of a connected $\k$-algebra,
$$\pi:T(x_1,\dots,x_N)\twoheadrightarrow A$$
be the projection. Then the following sequence of graded free left $A$-modules is exact:
$$A\cdot\{R_1,\dots,R_M\}\overset{d_2}\longrightarrow A\cdot\{X_1,\dots,X_N\}\overset{d_1}\longrightarrow A\overset{\varepsilon}\longrightarrow\k\to 0,$$
$$d_2(R_j):=-\sum_{i=1}^N \pi(\overline{r_{ji}})\cdot X_i,\quad d_1(X_i):=x_i.$$
\end{prp}
\begin{proof}
We first prove that the sequence is a chain complex. Indeed, $\varepsilon(d_1(X_i))=\varepsilon(x_i)=0$ and
$$d_1(d_2(R_j))=\sum_{i=1}^N\pi(r_{ji})\cdot d_1(X_i)=\sum_{i=1}^N\pi(r_{ji})x_i=\pi\left(\sum_{i=1}^N r_{ji}x_i\right)=\pi(r_j)=0\in A$$
($r_j\in\Ker\varepsilon,$ hence $r_j=\sum_i r_{ji}x_i$).
We check the exactness in the term $A.$ Let $y\in \Ker\varepsilon\subset A.$ We have $y=\pi(w)$ for some element $w\in T(x_1,\dots,x_N)$ of positive degree, hence
$$y=\pi\left(\sum_{i=1}^N w_ix_i\right)=\sum_{i=1}^N\pi(w_i)x_i=d_1\left(-\sum_{i=1}^N\pi(\overline{w}_i)\cdot X_i\right)\in\Img d_1.$$
Finally, we check the exactness in the term $A\cdot\{X_1,\dots,X_N\}.$ Suppose that $\sum_{i=1}^Na_i\cdot X_i\in\Ker d_1,$ so $\sum_{i=1}^N\overline{a}_ix_i=0.$ We have $a_i=\pi(v_i)$ for some $v_i\in T(x_1,\dots,x_N).$ Then the element $w:=\sum_{i=1}^N \overline{v}_ix_i\in T(x_1,\dots,x_N)$ belongs to $\Ker\pi.$ This kernel is a two-sided ideal generated by $r_j.$ Hence
$w=\sum_{j=1}^M\sum_\alpha u_{j,\alpha}r_j w_{j,\alpha}$
for some $u_{j,\alpha},w_{j,\alpha}\in T(x_1,\dots,x_N).$ We can rewrite it as
$$w=\sum_{j=1}^M\sum_\alpha u_{j,\alpha}r_j\varepsilon(w_{j,\alpha})+\sum_{j=1}^M\sum_\alpha\sum_{i=1}^N u_{j,\alpha}r_j w_{j,\alpha,i}x_i=\sum_{i=1}^N\sum_{j=1}^M\sum_\alpha \left(\varepsilon(w_{j,\alpha})u_{j,\alpha}r_{ji}+u_{j,\alpha}r_jw_{j,\alpha,i}\right)x_i.$$
On the other hand, $w=\sum_{i=1}^N\overline{v}_ix_i.$ Such representation is unique, so we have
$$\overline{v}_i=\sum_{j=1}^M\sum_\alpha \varepsilon(w_{j,\alpha})u_{j,\alpha}r_{ji}+u_{j,\alpha}r_jw_{j,\alpha,i},\quad i=1,\dots,N.$$
Applying $\pi$ to both parts of this identity, we obtain $\overline{a}_i=\sum_{j=1}^M\sum_\alpha\varepsilon(w_{j,\alpha})\pi(u_{j,\alpha})\pi(r_{ji}),$ since $\pi(v_i)=a_i$ and $\pi(r_j)=0.$
Finally,
\[
\sum_{i=1}^N a_i\cdot X_i=-\sum_{j=1}^M\sum_\alpha \varepsilon(w_{j,\alpha})\pi(\overline{u}_{j,\alpha})\pi(\overline{r}_{ji})\cdot X_i=d_2\left(-\sum_{j=1}^M\sum_\alpha\varepsilon(w_{j,\alpha}) \pi(u_{j,\alpha})\cdot R_j\right)\in\Img d_2.
\qedhere
\]
\end{proof}
\begin{rmk}
Proposition \ref{prp:presentation_exact_sequence} holds for presentations of \emph{augmented} algebras such that $\varepsilon(x_i)=\varepsilon(r_j)=0.$ The corresponding exact sequence is called the ``Koszul resolution'' in \cite[\S 2]{anick_dicks}.
\end{rmk}
\begin{crl}
\label{crl:generators_generate_Tor1}
Let $A=T(x_1,\dots,x_N)/(r_1,\dots,r_M)$ be a presentation of a connected graded $\k$-algebra, which is a free $\k$-module. Then the $\k$-module $\Tor^A_1(\k,\k)$ is additively generated by images of cycles $[x_1],\dots,[x_N]\in\oB_1(A).$
\end{crl}
\begin{proof}
We extend the exact sequence from Proposition~\ref{prp:presentation_exact_sequence} to a free resolution
$$\dots\to A\cdot\{X_1,\dots,X_N\}\overset{d_1}\longrightarrow A\overset{\varepsilon}\longrightarrow \k\to 0,\quad d_1(X_i)=x_i$$
of the left $A$-module $\k.$ Consider the diagram
$$\xymatrix{
\dots\ar[r]
&
A\cdot\{X_1,\dots,X_N\}
\ar[r]^-{d_1}
\ar[d]^-{X_i\mapsto [x_i]}
&
A
\ar[r]^-\varepsilon
\ar[d]^-{a\mapsto a[]}
&
\k
\ar[r]
\ar@{=}[d]
&
0\\
\dots
\ar[r]
&
\B_1(A)
\ar[r]^-{d_{\B,1}}
&
\B_0(A)
\ar[r]^-{\varepsilon}
&
\k
\ar[r]
&
0.
}$$
It is commutative, since $d_1(a\otimes X_i)=-\overline{a}x_i$ and $d_{\B,1}(a[x_i])=-\overline{a}x_i[].$ Hence it can be extended to a map of resolutions (e.g. using Lemma \ref{lmm:building_resolution}). Apply the functor $\k\otimes_A(-).$ We obtain a map of chain complexes
$$\xymatrix{
\dots\ar[d]\ar[r]
&
\k\cdot\{X_1,\dots,X_N\}
\ar[r]^-{0}
\ar[d]^-{X_i\mapsto [x_i]}
&
\k
\ar[r]
\ar[d]^-{a\mapsto a[]}
&
0\\
\dots
\ar[r]
&
\oB_1(A)
\ar[r]^-{d_{\oB,1}}
&
\oB_0(A)
\ar[r]
&
0.
}$$
The homology of both complexes equals $\Tor^A(\k,\k),$ and the induced map in homology is an isomorphism. The elements $X_i$ in the first complex are cycles, and their images generate $\Tor^A_1(\k,\k).$
\end{proof}
\begin{crl}
\label{crl:tensor_Tor}
Let $A=T(x_1,\dots,x_N)$ be the tensor algebra over a ring $\k$, where $x_1,\dots,x_N$ are homogeneous elements of positive degrees. Then $\Tor^A_1(\k,\k)$ is a free $\k$-module with the basis represented by cycles $[x_1],\dots,[x_N]\in\oB_1(A).$ Moreover, $\Tor^A_i(\k,\k)=0$ for $i>1.$
\end{crl}
\begin{proof}
By Proposition~\ref{prp:presentation_exact_sequence}, the sequence
$$0\to A\cdot\{X_1,\dots,X_N\}\overset{d_1}\longrightarrow A\overset{\varepsilon}\longrightarrow\k\to 0,\quad d_1(X_i)=x_i,$$
is exact. As in the proof of Corollary \ref{crl:generators_generate_Tor1}, we obtain a map of chain complexes
$$\xymatrix{
0\ar[d]\ar[r]
&
\k\cdot\{X_1,\dots,X_N\}
\ar[r]^-{0}
\ar[d]^-{X_i\mapsto [x_i]}
&
\k
\ar[r]
\ar[d]^-{a\mapsto a[]}
&
0\\
\dots
\ar[r]
&
\oB_1(A)
\ar[r]^-{d_{\oB,1}}
&
\oB_0(A)
\ar[r]
&
0.
}$$
Homology of both complexes is equal to $\Tor^A(\k,\k),$ and the induced map in homology is the identity.
\end{proof}

\subsection{A presentation that corresponds to cycles}

Recall that $\Tor^A(\k,\k)\cong H(\oB(A))$ is $A$ is a free $\k$-module.
The following lemma is proved by Lemaire \cite[Corollaire 1.2.3]{lemaire} in the case of field coefficients.
\begin{lmm}
\label{lmm:tor_bijectivity_criterion}
~\\
Let $f:A\to C$ be a morphism of connected $\k$-algebras, where $\k$ is a commutative ring with unit.
\begin{enumerate}
\item Suppose that the map $f_{*,1}:H_1(\oB(A))\to H_1(\oB(C))$ is surjective. Then $f:A\to C$ is surjective.
\item Suppose that $f_{*,1}:H_1(\oB(A))\to H_1(\oB(C))$ is bijective, and the map $f_{*,2}:H_2(\oB(A))\to H_2(\oB(C))$ is surjective. Then $f:A\to C$ is an isomorphism.
\end{enumerate}
\end{lmm}
We prove by induction that the maps $f_n:A_n\to C_n$ are surjective (bijective). The base case is the bijection $A_0\cong\k\cong C_0.$ Recall that the bar construction $\oB(A)$ is the chain complex
$$\dots\to \oB_3(A)\overset{d_{3}}\longrightarrow \oB_2(A)\overset{d_{2}}\longrightarrow \oB_1(A)\overset{0}\longrightarrow\k\to 0,$$
$$\oB_k(A)=I(A)^{\otimes k},\quad d_{2}(x\otimes y)=\overline{x}y,\quad d_{3}(x\otimes y\otimes z)=\overline{x}y\otimes z + \overline{x}\otimes \overline{y}z.$$
We denote $f_\#:\oB(A)\to\oB(C).$

\begin{proof}[Proof of statement (1)]
Suppose that $f:A_i\to C_i$ is surjective for $i<n.$ Consider the following map of exact sequences:
$$\xymatrix{
\oB_2(A)_n
\ar[r]^-{d_2}
\ar@{->>}[d]^-{f_{\#,2}}
&
\oB_1(A)_n\cong A_n
\ar[r]
\ar[d]^-{f}
&
H_1(\oB(A))_n
\ar[r]
\ar@{->>}[d]^-{f_{*,1}}
&
0
\ar@{=}[d]
\\
\oB_2(C)_n
\ar[r]^-{d_2}
&
\oB_1(C)_n\cong C_n
\ar[r]
&
H_1(\oB(C))_n
\ar[r]
&
0.
}$$
The map $f_{\#,2}$ is surjective, since it is a direct sum of maps $f\otimes f:A_i\otimes A_j\to C_i\otimes C_j$ for $i,j<n,$ and $f$ is surjective in these degrees by the inductive hypothesis. The surjectivity $f_{*,1}$ is given, and $0\to 0$ is injective. Hence $f:A_n\to C_n$ is surjective by the ``first half of five lemma'' \cite[Proposition 2.72(i)]{rotman}.
\end{proof}
\begin{proof}[Proof of statement (2)]
Suppose that $f:A_i\to C_i$ is bijective for $i<n.$ Consider the following map of exact sequences:
$$
\xymatrix{
\oB_3(A)_n
\ar[r]^-{d_{3}}
\ar@{->>}[d]^-{f_{\#,3}}
&
\Ker d_{2}
\ar[r]
\ar[d]^-{\varphi}
&
H_2(\oB(A))_n
\ar[r]
\ar@{->>}[d]^-{f_{*,2}}
&
0
\ar@{=}[d]
\\
\oB_3(C)_n
\ar[r]^-{d_3}
&
\Ker d_{2}
\ar[r]
&
H_2(\oB(A))_n
\ar[r]
&
0.
}$$
The map $f_{\#,3}$ is surjective, since it is a direct sum of maps $f\otimes f\otimes f:A_i\otimes A_j\otimes A_k\to C_i\otimes C_j\otimes C_k$ for $i,j,k<n,$ and $f$ is surjective in these degrees. The surjectivity of $f_{*,2}$ is given, and $0\to 0$ is injective. Hence $\varphi$ is surjective by the ``first half of five lemma''. Now consider the following map of exact sequences:
$$\xymatrix{
\Ker d_{\oB,2}
\ar[r]
\ar@{->>}[d]^-\varphi
&
\oB_2(A)_n
\ar[d]^-{f_{\#,2}}_-\simeq
\ar[r]^-{d_2}
&
\oB_1(A)_n\cong A_n
\ar[r]
\ar[d]^-{f}
&
H_1(\oB(A))_n
\ar@{^(->}[d]^-{f_{*,1}}
\\
\Ker d_{\oB,2}
\ar[r]
&
\oB_2(C)_n
\ar[r]^-{d_2}
&
\oB_1(C)_n\cong C_n
\ar[r]
&
H_1(\oB(C))_n.
}$$
We proved that $\varphi$ is surjective. The map $f_{\#,2}$ is bijective by the inductive hypothesis (it is a direct sum of maps $f\otimes f:A_i\otimes A_j\to C_i\otimes C_j,$ $i,j<n$); in particular, it is injective. The injectivity of $f_{*,1}$ is given. Hence the map $f:A_n\to C_n$ is injective by the ``second half of five lemma'' \cite[Proposition 2.72(ii)]{rotman}. By (1), this map is also surjective.
\end{proof}

The following theorem allows one to obtain a presentation of a connected $\k$-algebra $A$, knowing the structure of $\k$-modules $H_1(\oB(A))$ and $H_2(\oB(A)).$ In the proof, we do not use the notation $[x|y|z]$ for elements of the bar construction, and write $x\otimes y\otimes z$ instead. Therefore, $[c]$ always denotes the class in $H(\oB(\Gamma))$ represented by a cycle $c\in\oB(\Gamma).$

We also use the following notation. Let $a_1,\dots,a_N\in A$ be some homogeneous elements of positive degree and $K,L\in T(x_1,\dots,x_N)$ be homogeneous non-commutative polynomials that belong to the augmentation ideal. Then the elements $K(a_1,\dots,a_N),L(a_1,\dots,a_N)\in I(A)$ are defined, and hence we can consider the elements $K(a_1,\dots,a_N)\otimes L(a_1,\dots,a_N)\in I(A)\otimes I(A)=\oB_2(A)$ and
$$d_{\oB,2}(K(a_1,\dots,a_N)\otimes L(a_1,\dots,a_N))=\overline{K}(a_1,\dots,a_N)\cdot L(a_1,\dots,a_N)\in\oB_1(A)= I(A).$$
\begin{thm}
\label{thm:presentation_from_cycles}
Let $A$ be a connected algebra over a commutative ring $\k$ with unit.
\begin{enumerate}
\item Suppose that, for homogeneous elements $a_1,\dots,a_N\in A_{>0},$ the $\k$-module $H_1(\oB(A))$ is additively generated by the classes $[a_1],\dots,[a_N]\in H_1(\oB(A)).$ Then $A$ is multiplicatively generated by $a_1,\dots,a_N.$
\item Suppose that the $\k$-module $H_1(\oB(A))$ is additively generated by $N$ elements $[a_1],\dots,[a_N]$ modulo $R$ relations
$$\sum_{i=1}^N \lambda_{ri}[a_i]=0\in H_1(\oB(A)),\quad r=1,\dots,R,~\lambda_{ri}\in\k.$$
Suppose that homogeneous polynomials $P_{j,\alpha},Q_{j,\alpha},K_{r,\beta},L_{r,\beta}\in T(x_1,\dots,x_N)$ are such that
$$\sum_{i=1}^N\lambda_{ri}\cdot a_i=d_{\oB,2}\left(\sum_\beta K_{r,\beta}(a_1,\dots,a_N)\otimes L_{r,\beta}(a_1,\dots,a_N)\right)\in I(A),\quad r=1,\dots,R,$$
and the cycles in bar construction
$$\sum_\alpha P_{j,\alpha}(a_j,\dots,a_N)\otimes Q_{j,\alpha}(a_1,\dots,a_N)\in I(A)\otimes I(A),\quad j=1,\dots,M,$$
generate the $\k$-module $H_2(\oB(A)).$ Then the algebra $A$ has a presentation
$$A\cong T(x_1,\dots,x_N)/\left(\sum_{i=1}^N\lambda_{ri}x_i=\sum_\beta\overline{K}_{r,\beta}\cdot L_{r,\beta},~r=1,\dots,R;~\sum_\alpha\overline{P}_{j,\alpha}\cdot Q_{j,\alpha}=0,~j=1,\dots,M\right).$$
\end{enumerate}
(Here $N,M,R$ can be infinite of any cardinality.)
\end{thm}
\begin{proof}[Proof of statement (1)]
Consider the morphism
$f:T(x_1,\dots,x_N)\to A,$ $x_i\mapsto a_i,$ of connected algebras.
The classes $[a_1],\dots,[a_N]$ generate $H_1(\oB(A))$ and are images of classes $[x_1],\dots,[x_N]$ with respect to the map $f_{*,1}:H_1(\oB(T(x_1,\dots,x_N))\to H_1(\oB(A)).$ Hence $f_{*,1}$ is surjective. By Lemma \ref{lmm:tor_bijectivity_criterion}(1), $f$ is surjective.
\end{proof}
\begin{proof}[Proof of statement (2)]
Consider the algebra
$$C:=T(x_1,\dots,x_N)/\left(\sum_{i=1}^N\lambda_{ri}x_i=\sum_\beta\overline{K}_{r,\beta}\cdot L_{r,\beta},~r=1,\dots,R;~\sum_\alpha\overline{P}_{j,\alpha}\cdot Q_{j,\alpha}=0,~j=1,\dots,M\right).$$
The following identities in $A$ are given:
$$\sum_{i=1}^N \lambda_{ri}\cdot a_i=\sum_\beta\overline{K}_{r,\beta}(a_1,\dots,a_N)\cdot L_{r,\beta}(a_1,\dots,a_N),\quad 0=\sum_\alpha\overline{P}_{j,\alpha}(a_1,\dots,a_N)\cdot Q_{j,\alpha}(a_1,\dots,a_N).$$
Hence the morphism $f:C\to A,$ $x_i\mapsto a_i,$ is well defined. The induced map $f_{*,1}:H_1(\oB(C))\to H_1(\oB(A))$ is surjective, since the elements $[a_i]=f_{*,1}([x_i])$ generate $H_1(\oB(A)).$

We prove that $f_{*,1}$ is injective. Let $\xi\in H_1(\oB(C))$ and $f_{*,1}(\xi)=0.$ By Corollary \ref{crl:tensor_Tor} and surjectivity of $T(x_1,\dots,x_N)\to C$, we have $\xi=\sum_{i=1}^N\mu_i\cdot [x_i]$ for some $\mu_i\in\k.$ Then $0=f_*(\xi)=\sum_i\mu_i [a_i]\in H_1(\oB(A)).$ All linear relations between $[a_1],\dots,[a_N]$ follow from the relations $\sum_i\lambda_{ri}[a_i]=0,$ hence $\mu_i=\sum_{r=1}^R c_r\lambda_{ri}$ for some $c_r\in\k.$ It follows that $\xi$ is represented by the cycle
$$\sum_{i=1}^N\sum_{r=1}^R c_r\lambda_{ri}\cdot x_i=\sum_{r=1}^R c_r\sum_\beta\overline{K}_{r,\beta}\cdot L_{r,\beta}=d_{\oB,2}\left(\sum_{r=1}^R c_r\sum_\beta K_{r,\beta}\otimes L_{r,\beta}\right)\in\oB_1(C).$$
Hence $\xi=0.$ We proved that $f_{*,1}$ is bijective.

The elements $\sum_\alpha P_{i,\alpha}\otimes Q_{i,\alpha}\in I(C)\otimes I(C)$ are cycles in $\oB_2(C),$ and their images generate $H_2(\oB(A)).$ Hence $f_{*,2}:H_2(\oB(C))\to H_2(\oB(A))$ is surjective. Conditions of Lemma \ref{lmm:tor_bijectivity_criterion}(2) are satisfied, so $f$ is bijective.
\end{proof}

\subsection{Bounds on the number of homogeneous generators and relations}
Let $A$ be a connected $\k$-algebra. Proposition~\ref{prp:presentation_exact_sequence} gives a lower bound on the number of generators and relations in the homogeneous presentations of $A,$ and Theorem~\ref{thm:presentation_from_cycles} gives an upper bound. These bounds coincide if $\k$ is a principal ideal domain, $A$ is a free $\k$-module, and graded components are finitely generated. We introduce some notations.
\begin{dfn}
Let $M$ be a finitely generated module over a principal ideal domain $\k.$ By the structure theorem of such modules, we have
\begin{equation}
\label{eq:pid_module}
M\simeq \k/(d_1)\oplus\dots\oplus\k/(d_r),
\end{equation}
where $d_1,\dots,d_r\in\k$ are non-invertible, and $d_{i}\mid d_{i+1}$ for all $i=1,\dots,r-1.$ The number $r$ is determined uniquely, and the elements $d_i$ --- uniquely up to a multiplication by an invertible element. Hence, the numbers
$$\gen M:=r,\quad\rel M:=\max\{s:~d_s\neq 0\}$$
are well defined. We get a short exact sequence
$\k^{\rel M}\to\k^{\gen M}\to M\to 0.$
\end{dfn}
\begin{lmm}
\label{lmm:PID_smith}
Let $\k$ be a principal ideal domain. Suppose that there is a short exact sequence
$\k^A\overset{f}\longrightarrow\k^B\to M\to 0$ for some $A,B<\infty.$ Then $A\geq\rel M$ and $B\geq\gen M.$
\end{lmm}
\begin{proof}
We can assume that $f$ is in the Smith normal form, that is, $f$ is represented by a diagonal matrix with diagonal elements $d'_1,\dots,d'_s$ such that $d'_1\mid d'_2\mid\dots\mid d'_s.$ Remove all nonzero columns: this preserves cokernel and does not increase $A.$ If $d'_i$ is invertible, remove the $i$-th row and the $i$-th column: this preserves cokernel and diminish $A$ and $B$ by 1. We obtain a diagonal matrix $B'\times A'$ having no zero columns and no invertible elements on diagonal. Hence the cokernel is exactly of the form \eqref{eq:pid_module} for $B'=r=\gen M$ and $A'=s=\rel M.$ 
\end{proof}
\begin{lmm}
\label{lmm:PID_exact}
Let $\k$ be a principal ideal domain and $0\to \k^a\to\k^b\overset{f}\longrightarrow\k^c\to 0$ be an exact sequence of $\k$-modules for some $a,b,c<\infty.$ Then $b=a+c.$
\end{lmm}
\begin{proof}
We can assume that $f$ is in a Smith normal form. In this basis, $f$ is represented by a diagonal matrix $c\times b.$ Since $f$ is surjective, the matrix has no nonzero rows, and all diagonal elements are non-invertible. Hence $\Ker f\simeq\k^{b-c}.$ We have $\k^d\not\simeq\k^{d'}$ for $d\neq d',$ so $a=b-c.$
\end{proof}

Recall that we consider $G$-graded algebras that are connected with respect to the $\ZZ$-grading given by a map $G\to\ZZ.$ 

\begin{thm}
\label{thm:size_of_pres_PID}
Let $A$ be a connected associative algebra with unit over a principal ideal domain $\k.$ Suppose that $\k$-modules $\Tor^A_1(\k,\k)_n$ and $\Tor^A_2(\k,\k)_n$ are finitely generated for all $n\in G.$ Then
\begin{enumerate}
\item If $A$ is a free $\k$-module, it admits a homogeneous presentation that contains (for every $n$) precisely $\gen\Tor^A_1(\k,\k)_n$ generators and $\gen\Tor^A_2(\k,\k)_n+\rel\Tor^A_1(\k,\k)_n$ relations
of degree $n$.
\item If $A$ admits a homogeneous presentation that contains $N_n$ generators and $M_n$ relations of degree $n,$ then
\begin{equation}
\label{eq:gens_rels_lower_bound}
N_n\geq\gen\Tor^A_1(\k,\k)_n,
\quad
M_n\geq\gen\Tor^A_2(\k,\k)_n+
\rel\Tor^A_1(\k,\k)_n.
\end{equation}
\end{enumerate}
\end{thm}
\begin{proof}[Proof of statement (1)]
For every $n$, choose a set of $\gen(\Tor^A_1(\k,\k)_n)$ additive generators for the $\k$-module $\Tor^A_1(\k,\k)_n$, a set of $\rel(\Tor^A_1(\k,\k)_n)$ linear relations between them, and a set of $\gen(\Tor^A_2(\k,\k)_n)$ generators for $\Tor^A_2(\k,\k)_n.$ These elements are represented by cycles and boundaries in the bar construction. Applying Theorem~\ref{thm:presentation_from_cycles} to them, we obtain a presentation of required size.
\end{proof}
\begin{proof}[Proof of statement (2)]
Apply Proposition~\ref{prp:presentation_exact_sequence} and continue the exact sequence to the free resolution of the left $A$-module $\k.$ It has the form
$$\dots\to A\otimes\k^M\to A\otimes\k^N\to A\overset{\varepsilon}\longrightarrow\k\to 0.$$
Applying the functor $\k\otimes_A(-),$ we obtain a chain complex of graded $\k$-modules
$$\dots\to \k^M\overset{\partial}\longrightarrow \k^N\overset{0}\longrightarrow \k\to 0,$$
having $\Tor^A(\k,\k)$ as homology. Therefore, for some $\partial_n:\k^{M_n}\to\k^{N_n}$ we have
$$\Coker\partial_n\simeq \Tor^A_1(\k,\k)_n,\quad\Ker\partial_n\twoheadrightarrow\Tor^A_2(\k,\k)_n.$$
In particular, $\Tor^A_1(\k,\k)_n$ is generated by $N_n$ elements, so $N_n\geq\gen\Tor^A_1(\k,\k)_n.$

If $M_n$ is infinite, both inequalities \eqref{eq:gens_rels_lower_bound} are true, since the right side is finite. If $M_n$ is finite, then $N_n$ is finite, since $\Coker\partial_n$ is finitely generated. Thus $\Ker\partial_n\subset\k^{M_n},$ $\Img\partial_n\subset\k^{N_n}$ are submodules of finitely generated free modules, so these modules are free: $\Ker\partial_n\simeq\k^P,$ $\Img\partial_n\simeq\k^Q.$ We obtain exact sequences
$$\k^P\to\k^{N_n}\to\Tor^A_1(\k,\k)_n\to 0,\quad \k^{Q}\to\Tor^A_2(\k,\k)_n\to 0,\quad 0\to \k^P\to \k^{M_n}\to \k^Q\to 0.$$
Then $N_n\geq\gen\Tor^A_1(\k,\k)_n,$ $P\geq\rel\Tor^A_1(\k,\k)_n,$ $Q\geq\gen\Tor^A_2(\k,\k)_n$ by Lemma \ref{lmm:PID_smith} and $P+Q=M_n$ by Lemma \ref{lmm:PID_exact}. This proves the inequalities \eqref{eq:gens_rels_lower_bound}.
\end{proof}
As a corollary, we obtain a well known result by Wall \cite[\S 7]{wall}:
\begin{crl}
Let $A$ be a connected associative algebra with unit over a field $\k$. Then
\begin{enumerate}
\item $A$ admits a homogeneous presentation that contains (for every n) precisely $\dim_\k\Tor^A_1(\k,\k)_n$ generators and $\dim_\k\Tor^A_2(\k,\k)_n$ relations of degree $n.$
\item If $A$ admits a homogeneous presentation that contains $N_n$ generators and $M_n$ relations of degree $n,$ then
$N_n\geq\dim_\k\Tor^A_1(\k,\k)_n$ and $M_n\geq\dim_\k\Tor^A_2(\k,\k)_n.$\qed
\end{enumerate}
\end{crl}
We also obtain a criterion of freeness.
\begin{crl}[{\cite[Proposition 8.5.4]{neisendorfer_book}}]
\label{crl:PID_freeness_criterion}
Let $A$ be a connected associative algebra with unit over a principal ideal domain $\k$, which is a free $\k$-module. The following conditions are equivalent.
\begin{enumerate}
\item[(a)] $A$ is a free algebra (a tensor algebra on homogeneous generators).
\item[(b)] $\Tor^A_1(\k,\k)$ is a free $\k$-module, and $\Tor^A_2(\k,\k)=0.$
\end{enumerate}
\end{crl}
\begin{proof}
By Corollary~\ref{crl:tensor_Tor}, (a) implies (b). Conversely, suppose that (b) holds. Then $\rel\Tor^A_1(\k,\k)=\gen\Tor^A_2(\k,\k)=0.$ By Theorem~\ref{thm:size_of_pres_PID}(1), the algebra $A$ admits a presentation with no relations. Hence $A$ is free.
\end{proof}
\section{Loop homology and extensions of Hopf algebras}
\label{sec:appendix_section}
Consider a homotopy fibration $F\to E\overset{p}\longrightarrow B$ of simply connected spaces, such that $\Omega p:\Omega E\to \Omega B$ admits a homotopy section (i.e. there is a continuous map $\sigma:\Omega B\to \Omega E$ that preserves basepoints, and a homotopy $\Omega p\circ \sigma\sim\id_{\Omega B}$). It is well known that then $\Omega E$ is homotopy equivalent to $\Omega F\times\Omega B$ (see \cite[Theorem 5.2]{eckmann_hilton} and \cite[Proposition A.2]{brace_product}). If $\k$-homology of these loop spaces is free, we obtain an extension of Hopf algebras $\k\to H_*(\Omega F;\k)\to H_*(\Omega E;\k)\to H_*(\Omega B;\k)\to\k.$ In Theorem \ref{thm:section_hopf_extension} we give a full proof of this folklore result. We consider ordinary loop spaces instead of Moore loop spaces, so that $\Omega(X\times Y)\cong\Omega X\times\Omega Y$ is a strict isomorphism of H-spaces.

We have a natural isomorphism $$\alpha:\pi_n(A\times B)\overset\cong\longrightarrow\pi_n(A)\oplus\pi_n(B),\quad [f]\mapsto [\mathrm{pr}_A\circ f]\oplus [\mathrm{pr}_B\circ f]$$ for any $A,B$ and $n\geq 1.$ We denote basepoint inclusion by $\varepsilon:\ast\to Y$ and collapse map by $\eta:Y\to\ast.$
\begin{lmm}
\label{lmm:product_pi_map}
Let $X$ be a simply connected space and $\mu:\Omega X\times\Omega X\to\Omega X$ be the composition of loops. Then the following diagram is commutative:
$$\xymatrix{
\pi_n(\Omega X\times\Omega X)\ar[r]^-{\mu_*}\ar[d]_\alpha^-\cong & \pi_n(\Omega X)\\
\pi_n(\Omega X)\oplus\pi_n(\Omega X)\ar[ur]_-{~~(x,y)\mapsto x+y}
}$$
\end{lmm}
\begin{proof}
Let elements $x,y\in\pi_n(\Omega X$) be represented by maps $f,g:S^n\to\Omega X.$ Consider the element $z=[f\times\eta]+[\eta\times g]\in\pi_n(\Omega X\times\Omega X).$

The map $\mu\circ(f\times\eta)$ is the composition $S^n\overset{f}\longrightarrow \Omega X\overset{\id\times\eta}\longrightarrow\Omega X\times\Omega X\overset{\mu}\longrightarrow\Omega X.$ The composition of two right maps is homotopic to the identity, hence $\mu\circ(f\times\eta)\sim f.$ Passing to homotopy groups, we have $\mu_*([f\times\eta])=x.$ Similarly, $\mu_*([\eta\times g])=y,$ hence $\mu_*(z)=x+y.$ On the other hand, $\alpha([f\times\eta])=[\mathrm{pr}_1\circ(f\times\eta)]\oplus [\mathrm{pr}_2\circ(f\times\eta)]=[f]\oplus[\eta\varepsilon]=x\oplus 0.$ Similarly, $\alpha([\eta\times g])=0\oplus y,$ hence $\alpha(z)=x\oplus y.$ We obtained $\mu_*(\alpha^{-1}(x\oplus y))=\mu_*(z)=x+y,$ so the diagram commutes.
\end{proof}

In the following lemma, we say that diagram commutes if it homotopy commutes.
\begin{lmm}
\label{lmm:section_homotopy_equivalence}
Let $F\overset{i}\longrightarrow E\overset{p}\longrightarrow B$ be a fibration of simply connected spaces, and $\sigma:\Omega B\to\Omega E$ be a homotopy section for $\Omega p.$
Consider the composition
$$f:\Omega F\times\Omega B\overset{\Omega i\times\sigma}\longrightarrow \Omega E\times \Omega E\overset{\mu}\longrightarrow \Omega E.$$
Then
\begin{enumerate}
\item $f$ is a weak homotopy equivalence;
\item $f$ respects the inclusion and the projection, that is, the following diagram commutes:
$$\xymatrix{
&
\Omega E
\ar[dr]^-{\Omega p}
&
\\
\Omega F\times\ast
\ar[ur]^-{\Omega i}
\ar[r]_-{\id\times\eta}
&
\Omega F\times\Omega B
\ar[u]^-f
\ar[r]_-{\varepsilon\times\id}
&
\ast\times\Omega B;
}$$
\item $f$ respects the left action of $\Omega F,$ that is, the following diagram commutes:
$$\xymatrix{
\Omega F\times\Omega F\times\Omega B
\ar[r]^-{\Omega i\times f}
\ar[d]^-{\mu\times\id}
&
\Omega E\times\Omega E
\ar[d]^-\mu
\\
\Omega F\times\Omega B
\ar[r]^-f
&
\Omega E;
}$$
\item $f$ respects the right coaction of $\Omega B$, that is, the following diagram commutes:
$$\xymatrix{
\Omega F\times\Omega B
\ar[r]^-{f}
\ar[dd]^-{\id\times\Delta}
&
\Omega E
\ar[d]^-{\Delta}
\\
&
\Omega E\times\Omega E
\ar[d]^-{\id\times\Omega p}
\\
\Omega F\times\Omega B\times\Omega B
\ar[r]^-{f\times\id}
&
\Omega E\times\Omega B.
}$$
\end{enumerate}
\end{lmm}
\begin{proof}
We have an exact sequence
$$\dots\to \pi_n(\Omega F)\overset{(\Omega i)_*}\longrightarrow\pi_n(\Omega E)\overset{(\Omega p)_*}\longrightarrow\pi_n(\Omega B)\to\dots,$$
where the map $(\Omega p)_*$ has a section $\sigma_*.$ For every $n\geq 1,$ we obtain an isomorphism of groups
$$\varphi:\pi_n(\Omega F)\oplus\pi_n(\Omega B)\overset{\simeq}\longrightarrow\pi_n(\Omega E),\quad \varphi(x,y)=(\Omega i)_*(x)+\sigma_*(y).$$
(We use that $\pi_1(\Omega X)$ is abelian.) By Lemma~\ref{lmm:product_pi_map} and naturality of $\alpha:\pi_n(\Omega F\times\Omega B)\to\pi_n(\Omega F)\oplus\pi_n(\Omega B)$ we have
$$\varphi\circ\alpha=(\mu\circ (\Omega i\times \sigma))_*=f_*:\pi_n(\Omega F\times\Omega B)\to \pi_n(\Omega E).$$
Hence $f_*$ is an isomorphism for all $n,$ so $f$ is a weak homotopy equivalence. Now consider the diagram
$$\xymatrix{
&
\Omega E
\ar[r]^-{\Omega p}
&
\Omega B
\\
\Omega E\times\ast
\ar[r]^-{\id\times\eta}
\ar[ur]^-{\id}
&
\Omega E\times\Omega E
\ar[u]^-\mu
\ar[r]^-{\Omega p\times\Omega p}
&
\Omega B\times\Omega B
\ar[u]^-{\mu}
\\
\Omega F\times\ast
\ar[u]^-{\Omega i\times\id}
\ar[r]^-{\id\times\eta}
&
\Omega F\times\Omega B
\ar[u]^-{\Omega i\times\sigma}
\ar[r]^-{\varepsilon\times\id}
&
\ast\times\Omega B
\ar[u]^-{\eta\times\id}
}$$
The triangle commutes, since $\eta$ is a homotopy unit in $\Omega E.$ The upper right square commutes, since $\Omega p$ is a map of H-spaces. The bottom left square commutes, since $\eta_{\Omega E}=\sigma\circ\eta_{\Omega B}:\ast\to\Omega E.$ Finally, the commutativity of bottom right square is equivalent to the existence of homotopies $\Omega p\circ\Omega i\sim\eta\circ\varepsilon$ and $\Omega p\circ\sigma\sim\id.$
The first homotopy exists, since $p\circ i$ is homotopy trivial; the second exists, since $\sigma$ is a homotopy section for $\Omega p.$ Hence the whole diagram is commutative.  The right side of the diagram is homotopic to $\id:\Omega B\to\Omega B,$ since $\eta$ is a homotopy unit in $\Omega B.$ We obtain a commutative diagram
$$\xymatrix{
&
\Omega E
\ar[r]^-{\Omega p}
&
\Omega B\\
\Omega F\times\ast
\ar[ur]^-{\Omega i\times\id}
\ar[r]^-{\id\times\eta}
&
\Omega F\times\Omega B
\ar[u]_-f
\ar[r]^-{\varepsilon\times\id}
&
\ast\times\Omega B
\ar[u]^-{\id}
}$$
that is equivalent to the diagram from (2).
Now consider the diagram
$$\xymatrix{
\Omega F\times\Omega F\times\Omega B
\ar[rr]^-{\Omega i\times\Omega i\times\sigma}
\ar[d]^-{\mu\times\id}
&&
\Omega E\times\Omega E\times\Omega E
\ar[r]^-{\id\times\mu}
\ar[d]^-{\mu\times\id}
&
\Omega E\times\Omega E
\ar[d]^-\mu
\\
\Omega F\times\Omega B
\ar[rr]^-{\Omega i\times\sigma}
&&
\Omega E\times\Omega E
\ar[r]^-\mu
&
\Omega E.
}$$
The left square commutes, since $\Omega i:\Omega F\to\Omega E$ is a map of H-spaces; the right square commutes, since $\mu$ is homotopy associative. The top side of the diagram equals $\Omega i\times (\mu\circ (\Omega i\times\sigma))=\Omega i\times f,$
the bottom side equals $f.$ Hence, it is the diagram from (3). Finally, consider the diagram
$$\xymatrix{
\Omega F\times\Omega B
\ar[rr]^-{\Omega i\times\sigma}
\ar[d]^-D
&&
\Omega E\times\Omega E
\ar[r]^-\mu 
\ar[d]^-D
&
\Omega E
\ar[d]^-\Delta
\\
\Omega F\times\Omega B\times\Omega F\times\Omega B
\ar[d]^-{\mathrm{pr}_{124}}
\ar[drr]^-{\Omega i\times\sigma\times\varepsilon\times\id}
\ar[rr]^-{\Omega i\times\sigma\times\Omega i\times\sigma}
&&
\Omega E\times\Omega E\times\Omega E\times\Omega E
\ar[r]^-{\mu\times\mu}
\ar[d]^-{\id\times\id\times\Omega p\times\Omega p}
&
\Omega E\times\Omega E
\ar[d]^-{\id\times\Omega p}
\\
\Omega F\times\Omega B\times\Omega B
\ar[rr]^-\phi
&&
\Omega E\times\Omega E\times\Omega B\times\Omega B
\ar[r]^-{\mu\times\mu}
&
\Omega E\times\Omega B,
}$$
where $\Delta(x):=(x,x),$ $D(x,y):=(x,y,x,y)$
and $\phi(f,b_1,b_2):=(\Omega i(f),\sigma(b_1),\ast,b_2).$
Clearly, the top two squares commute. The bottom right square commutes, since $\id:\Omega E\to\Omega E$ and $\Omega p:\Omega E\to\Omega B$ are maps of H-spaces. The upper triangle commutes, since $\Omega p\circ\Omega i\sim\varepsilon$ and $\Omega p\circ\sigma\sim\id;$ the bottom triangle commutes by the definition of $\phi$. The outer maps in the diagram give the required diagram (4).
\end{proof}

In the proof of next theorem we use the K\"unneth map $\kappa:H_*(X;\k)\otimes H_*(Y;\k)\to H_*(X\times Y;\k).$ It is natural and associative. It is an isomorphism if $H_*(Y;\k)$ is a free $\k$-module.

If $X$ is a simply connected space and $H_*(\Omega X;\k)$ is free over $\k$, this module is a connected $\k$-Hopf algebra with the standard cup coproduct (see Subsection \ref{subsec:hopf_extension}) and the Pontryagin product
$$m:H_*(\Omega X;\k)\otimes H_*(\Omega X;\k)\overset{\kappa}\longrightarrow H_*(\Omega X\times\Omega X;\k)\overset{\mu_*}\longrightarrow H_*(\Omega X;\k),$$ The unit and counit $\k\overset{\eta}\longrightarrow H_*(\Omega X;\k)\overset{\varepsilon}\longrightarrow \k$ are induced by the H-space maps $\ast\overset{\eta}\longrightarrow \Omega X\overset{\varepsilon}\longrightarrow\ast.$

\begin{thm}
\label{thm:section_hopf_extension}
Let $\k$ be an associative ring with unit.
Let
$F\overset{i}\longrightarrow E\overset{p}\longrightarrow B$ be a homotopy fibration of simply connected spaces such that $H_*(\Omega B;\k)$ and $H_*(\Omega F;\k)$ are free $\k$-modules, and the map $\Omega p$ admits a homotopy section $\sigma:\Omega B\to \Omega E.$
Consider the composition
$$\Phi:H_*(\Omega F;\k)\otimes H_*(\Omega B;\k)\overset{(\Omega i)_*\otimes\sigma_*}\longrightarrow
H_*(\Omega E;\k)\otimes H_*(\Omega E;\k)\overset{m}\longrightarrow
H_*(\Omega E;\k).
$$
Then
\begin{enumerate}
\item $\Phi$ is an isomorphism of $\k$-modules;
\item $(\Omega i)_* = \Phi\circ(\id_{H_*(\Omega F;\k)}\otimes\eta_{H_*(\Omega B;\k)});$
\item $(\Omega p)_*\circ\Phi = \varepsilon_{H_*(\Omega F;\k)}\otimes\id_{H_*(\Omega B;\k)};$
\item $\Phi$ is a morphism of left $H_*(\Omega F;\k)$-modules and right $H_*(\Omega B;\k)$-comodules, where the (co)module structure on $H_*(\Omega E;\k)$ is induced by the maps $(\Omega i)_*$ and $(\Omega p)_*.$
\end{enumerate}
In particular, $\k\to H_*(\Omega F;\k)
\overset{(\Omega i)_*}\longrightarrow
H_*(\Omega E;\k)
\overset{(\Omega p)_*}\longrightarrow
H_*(\Omega B;\k)\to\k
$ is an extension of connected Hopf algebras over $\k.$
\end{thm}
\begin{proof}
We write $H_*(\Omega X)$ instead of $H_*(\Omega X;\k).$ Note that $\sigma$ is continuous, and $\Omega i,$ $\Omega p$ are maps of H-spaces. Hence $\sigma_*$ is a map of coalgebras, and $(\Omega i)_*,$ $(\Omega p)_*$ are maps of Hopf algebras. By the naturality of K\"unneth map, the following diagram commutes:
$$\xymatrix{
H_*(\Omega F)\otimes H_*(\Omega B)
\ar[r]^-{(\Omega i)_*\otimes\sigma_*}
\ar[d]_-\kappa^-\simeq
&
H_*(\Omega E)\otimes H_*(\Omega E)
\ar[r]^-{m}
\ar[d]_-\kappa
&
H_*(\Omega E)
\ar@{=}[d]
\\
H_*(\Omega F\times\Omega B)
\ar[r]^-{(\Omega i\times\sigma)_*}
&
H_*(\Omega E\times\Omega E)\ar[r]^-{\mu_*}
&
H_*(\Omega E).
}$$
The top side of diagram equals $\Phi,$ the bottom side equals $f_*.$ Hence $\Phi$ is the composition $H_*(\Omega F)\otimes H_*(\Omega B)\overset{\kappa}\longrightarrow H_*(\Omega F\times\Omega B)\overset{f_*}\longrightarrow H_*(\Omega E).$
The left map is bijective by the assumption, the right map is bijective by Lemma \ref{lmm:section_homotopy_equivalence}(1). Hence $\Phi$ is an isomorphism, so (1) is proved. Consider the diagram
$$\xymatrix{
&
H_*(\Omega E)
\ar@/^/[dr]^-{(\Omega p)_*}
&
\\
H_*(\Omega F\times\ast)
\ar@/^/[ur]^-{(\Omega i)_*}
\ar[r]_-{(\id\times\eta)_*}
&
H_*(\Omega F\times\Omega B)
\ar[u]^-{f_*}
\ar[r]_-{(\varepsilon\times\id)_*}
&
H_*(\ast\times\Omega B)
\\
H_*(\Omega F)\otimes\k
\ar[u]^-{\kappa}_{\simeq}
\ar[r]^-{\id\otimes\eta}
&
H_*(\Omega F)\otimes H_*(\Omega B)
\ar[u]^-{\kappa}
\ar[r]^-{\varepsilon\otimes\id}
&
\k\otimes H_*(\Omega B).
\ar[u]^-{\kappa}_{\simeq}
}$$
The top half of the diagram commutes by Lemma \ref{lmm:section_homotopy_equivalence}(2), the bottom half commutes by naturality of $\kappa.$ Since $f_*\circ\kappa=\Phi,$ we have a commutative diagram
$$\xymatrix{
&
H_*(\Omega E)
\ar@/^/[dr]^-{(\Omega p)_*}
&\\
H_*(\Omega F)
\ar@/^/[ur]^-{(\Omega i)_*}
\ar[r]^-{\id\otimes\eta}
&
H_*(\Omega F)\otimes H_*(\Omega B)
\ar[u]^-{\Phi}
\ar[r]^-{\varepsilon\otimes\id}
&
H_*(\Omega B),
}$$
which proves (2) and (3). Now consider the diagram
$$\xymatrix{
H_*(\Omega F)\otimes H_*(\Omega F)\otimes H_*(\Omega B)
\ar[d]^-{\kappa\otimes\id}
\ar[r]^-{\id\otimes\kappa}
&
H_*(\Omega F)\otimes H_*(\Omega F\times\Omega B)
\ar[d]^-\kappa
\ar[r]^-{(\Omega i)_*\otimes f_*}
&
H_*(\Omega E)\otimes H_*(\Omega E)
\ar[d]^-\kappa
\\
H_*(\Omega F\times\Omega F)\otimes H_*(\Omega B)
\ar[r]^-{\kappa\otimes\id}
\ar[d]^-{\mu_*\otimes\id}
&
H_*(\Omega F\times\Omega F\times\Omega B)
\ar[r]^-{(\Omega i\times f)_*}
\ar[d]^-{(\mu\times\id)_*}
&
H_*(\Omega E\times\Omega E)
\ar[d]^-{\mu_*}
\\
H_*(\Omega F)\otimes H_*(\Omega B)
\ar[r]^-\kappa
&
H_*(\Omega F\times\Omega B)
\ar[r]^-{f_*}
&
H_*(\Omega E).
}$$
The bottom right square commutes by Lemma \ref{lmm:section_homotopy_equivalence}(3), the other squares commute by naturality of $\kappa.$ Since $\mu_*\circ\kappa=m:H_*(\Omega X)\otimes H_*(\Omega X)\to H_*(\Omega X),$ the outer maps in the diagram are
$$\xymatrix{
H_*(\Omega F)\otimes H_*(\Omega F)\otimes H_*(\Omega B)
\ar[rr]^-{(\Omega i)_*\otimes\Phi}
\ar[d]^-{m\otimes\id}
&&
H_*(\Omega E)\otimes H_*(\Omega E)
\ar[d]^-m
\\
H_*(\Omega F)\otimes H_*(\Omega B)
\ar[rr]^-\Phi
&&
H_*(\Omega E).
}$$
Hence $\Phi$ is a map of left $H_*(\Omega F)$-modules. Similarly, by Lemma \ref{lmm:section_homotopy_equivalence}(4) and the K\"unneth isomorphisms we have the commutative diagram
$$\xymatrix{
H_*(\Omega F)\otimes H_*(\Omega B)
\ar[r]^-{\Phi}
\ar[dd]^-{\id\times\Delta}
&
H_*(\Omega E)
\ar[d]^-{\Delta}
\\
&
H_*(\Omega E)\otimes H_*(\Omega E)
\ar[d]^-{\id\otimes(\Omega p)_*}
\\
H_*(\Omega F)\otimes H_*(\Omega B)\otimes H_*(\Omega B)
\ar[r]^-{\Phi\otimes\id}
&
H_*(\Omega E)\otimes H_*(\Omega B),
}$$
hence $\Phi$ is a map of right $H_*(\Omega B$)-comodules.

Since (1)-(4) hold, the maps of Hopf algebras $(\Omega i)_*:H_*(\Omega F)\to H_*(\Omega E)$ and $(\Omega p)_*:H_*(\Omega E)\to H_*(\Omega B)$ form an extension of Hopf algebras by Proposition \ref{prp:hopf_extension_criterion}. 
\end{proof}

Recall that an element $x\in A$ of a Hopf algebra is \emph{primitive} if $\Delta x = 1\otimes x + x\otimes 1.$ The set of primitive elements is a Lie subalgebra $PA\subset A.$ Every map of Hopf algebras $f:A\to A'$ induces a map of Lie algebras $Pf:=f|_{PA}:PA\to PA'.$
\begin{crl}
\label{crl:primitive_in_kernel}
Suppose that the conditions of Theorem~\ref{thm:section_hopf_extension} are met. Let $x\in H_*(\Omega E;\k)$ be a primitive element such that $(\Omega p)_*(x)=0.$ Then $x=(\Omega i)_*(y)$ for some $y\in H_*(\Omega F;\k).$
\end{crl}
\begin{proof}
Since $\k\to H_*(\Omega F;\k)\to H_*(\Omega E;\k)\to H_*(\Omega B;\k)\to\k$ is a Hopf algebra extension, the sequence $0\to PH_*(\Omega F;\k)\to PH_*(\Omega E;\k)\to PH_*(\Omega B;\k)$
is exact, see \cite[Proposition 4.10]{milnor_moore}. (This also easily follows from definitions).
We have $x\in\Ker(PH_*(\Omega E;\k)\to PH_*(\Omega B;\k))=\Img(PH_*(\Omega F;\k)\to PH_*(\Omega E;\k)).$
\end{proof}

\section{Commutator identities}
\label{sec:appendix_identities}
Fix elements $u_1,\dots,u_m$ of degree 1 in a graded associative algebra $\Gamma.$ For a subset $I=\{i_1<\dots<i_k\}\subset[m],$ we denote
$$\widehat{u}_I:= u_{i_1}\cdot\dotso\cdot u_{i_k},\quad c(I,x):=[u_{i_1},[u_{i_2},[\dots [u_{i_k},x]\dots]]],~x\in\Gamma.$$
We write $A<B$ when $A,B\subset[m]$ and $\max(A)<\min(B).$ If $A<B$, we have $\widehat{u}_{A\sqcup B}=\widehat{u}_A\cdot\widehat{u}_B$ and $c(A\sqcup B,x)=c(A,c(B,x)).$ Also, $\widehat{u}_\varnothing=1,$ $c(\varnothing,x)=x.$

Define the Koszul sign by $\theta(A,B):=|\{(a,b)\in A\times B:~a>b\}|$. In a graded commutative algebra, we would have $\widehat{u}_A\cdot\widehat{u}_B=(-1)^{\theta(A,B)}\widehat{u}_{A\sqcup B}$ if $A\cap B=\varnothing.$ It has the following properties:
\begin{enumerate}
\item $\theta(A,B)\equiv |A|\cdot|B|+\theta(B,A)\mod 2;$
\item If $A_1\sqcup B_1 < A_2\sqcup B_2,$ then
$$\theta(A_1\sqcup A_2,B_1\sqcup B_2)\equiv \theta(A_1,B_1)+\theta(A_2,B_2)+|A_2|\cdot |B_1|.$$
\end{enumerate}

For $I\subset[m],$ $j\in [m],$ we write $I_{<j}=\{i\in I:~i<j\},$ $I_{>j}=\{i\in I:~i>j\}.$ We also use $i$ as a shortened notation for $\{i\}.$

\subsection{Regrouping of monomials}
The following formulas can be used to express any monomial on $u_1,\dots,u_m$ as a linear combination of $c_{1}\cdot\dotso\cdot c_{s}\cdot\widehat{u}_B,$ $c_i=c(A_i,u_{j_i}),$ $A_i\neq\varnothing.$
\begin{lmm}
Let $I\subset[m],$ and let $x\in\Gamma$ be homogeneous. Then
\begin{equation}
\label{eq:u_Ix}
\widehat{u}_I\cdot x = \sum_{I=A\sqcup B}(-1)^{\theta(A,B)+\deg(x)\cdot |B|}c(A,x)\widehat{u}_B.
\end{equation}
\end{lmm}
\begin{proof}
Denote $d:=\deg(x).$ Induction on $|I|.$ The base $I=\varnothing$ is clear. The inductive step: let $i=\min(I),$ $I'=I\setminus i.$
Then the right hand side is equal to
\begin{multline*}
\sum_{I'=A\sqcup B}
(-1)^{\theta(i\sqcup A,B)+d\cdot |B|}
c(i\sqcup A,x)\widehat{u}_B
+
\sum_{I'=A\sqcup B}
(-1)^{\theta(A,i\sqcup B)+d\cdot|i\sqcup B|}
c(A,x)\widehat{u}_{i\sqcup B}
\\=
\sum_{I'=A\sqcup B}
(-1)^{\theta(A,B)+d\cdot |B|}\left(
[u_i,c(A,x)]+(-1)^{|A|+d}c(A,x)u_i
\right)\cdot\widehat{u}_B
\\=
\sum_{I'=A\sqcup B}
(-1)^{\theta(A,B)+d\cdot|B|}\,
u_i c(A,x)\cdot\widehat{u}_B.
\end{multline*}
By the inductive hypothesis, this sum is equal to $u_i\cdot \widehat{u}_{I'}x=\widehat{u}_I\cdot x.$
\end{proof}
\begin{prp}
\label{prp:u_Iu_j}
Let $I\subset[m],$ $j\in [m].$ Then
$$\widehat{u}_I\cdot u_j =\sum_{\begin{smallmatrix}
I=A\sqcup B:\\
\max(A)>j
\end{smallmatrix}}
(-1)^{\theta(A,B)+|B|}c(A,u_j)\widehat{u}_B+
(-1)^{|I_{>j}|}\cdot\begin{cases}
\widehat{u}_{I\sqcup j},& j\notin I;\\
\widehat{u}_{I_{<j}}\cdot u_j^2\cdot\widehat{u}_{I_{>j}},& j\in I.
\end{cases}
$$
\end{prp}
\begin{proof}
Denote $P=I_{\leq j},$ $Q=I_{>j}.$ Then $P<Q,$ therefore
$$
\widehat{u}_I=\widehat{u}_P\widehat{u}_Q,
\quad
r:=
\widehat{u}_P\,
u_j\,
\widehat{u}_Q
=
\begin{cases}
\widehat{u}_{I\sqcup\{j\}},
&j\notin I;\\
\widehat{u}_{I_{<j}}\cdot
u_j^2\cdot
\widehat{u}_{I_{>j}},
&j\in I.
\end{cases}$$
Apply the formula \eqref{eq:u_Ix} to $\widehat{u}_Q\cdot u_j,$ and consider the summand with $A_2=\varnothing$ separately:
\begin{multline*}
\widehat{u}_I\cdot u_j=
\widehat{u}_P\,
\widehat{u}_Q\,
u_j
=
\sum_{Q=A_2\sqcup B_2}
(-1)^{\theta(A_2,B_2)+|B_2|}\,
\widehat{u}_P\,
c(A_2,u_j)\,
\widehat{u}_{B_2}
\\=
(-1)^{|Q|}
\widehat{u}_P\,
u_j\,
\widehat{u}_Q
+
\sum_{\begin{smallmatrix}
Q=A_2\sqcup B_2:\\
A_2\neq\varnothing
\end{smallmatrix}}
(-1)^{\theta(A_2,B_2)+|B_2|}\,
\widehat{u}_P\,
c(A_2,u_j)\,
\widehat{u}_{B_2}.
\end{multline*}
Applying \eqref{eq:u_Ix} to $\widehat{u}_{P}\cdot c(A_2,u_j)$, we obtain the required identity:
\begin{multline*}
\widehat{u}_I\cdot u_j=
(-1)^{|Q|}r
+
\sum_{P=A_1\sqcup B_1}
\sum_{\begin{smallmatrix}
Q=A_2\sqcup B_2:\\
A_2\neq\varnothing
\end{smallmatrix}}
(-1)^{\theta(A_1,B_1)+(|A_2|+1)\cdot |B_1|+
\theta(A_2,B_2)+|B_2|}
\,c(A_1,c( A_2,u_j))\widehat{u}_{B_1}\widehat{u}_{B_2}
\\=
(-1)^{|Q|}r+\sum_{\begin{smallmatrix}
P\sqcup Q = A\sqcup B:\\
A_{>j}\neq\varnothing
\end{smallmatrix}}
(-1)^{\theta(A,B)+|B|}c(A,u_j)\widehat{u}_B.
\qedhere
\end{multline*}
\end{proof}

\subsection{Identities for nested commutators}
In this section $\Gamma$ can be a Lie superalgebra.
\begin{lmm}
For $I\subset[m]$ and homogeneous elements $x,y\in\Gamma,$ we have
\begin{multline}
\label{c(I,[x,y])_expanded}
c(I,[x,y])=
\sum_{I=A\sqcup B}
(-1)^{\theta(A,B)+\deg(x)\cdot |B|}
\left[c(A,x),c(B,y)\right]
\\=
[c(I,x),y]+
(-1)^{\deg(x)\cdot|I|}
[x,c(I,y)]+
\sum_{\begin{smallmatrix}
I=A\sqcup B,\\
A,B\neq\varnothing
\end{smallmatrix}}
(-1)^{\theta(A,B)+\deg(x)\cdot |B|}
\left[c(A,x),c(B,y)\right].
\end{multline}
\end{lmm}
\begin{proof}
The second identity follows from $\theta(\varnothing,I)=\theta(I,\varnothing)=0$ and $c(\varnothing,x)=x.$ Let us prove the first identity by induction on $|I|.$ The base $I=\varnothing$ is clear. The inductive step: denote $i=\min(I),$ $I'=I\setminus i,$ $d=\deg(x).$ Then, by the inductive hypothesis,
\begin{multline*}
c(I,[x,y])=[u_i,c(I',[x,y])]=
\sum_{I'=A'\sqcup B'}
(-1)^{\theta(A',B')+d\cdot |B'|}
[u_i,[c(A',x),c(B',y)]]
\\=
\sum_{I'=A'\sqcup B'}
(-1)^{\theta(A',B')+d\cdot|B'|}
[[u_i,c(A',x)],c(B',y)]
+
\sum_{I'=A'\sqcup B'}
(-1)^{\theta(A',B')+d\cdot|B'|+d+|A'|}[c(A',x),[u_i,c(B',y)]]
\\=
\sum_{I'=A'\sqcup B'}
(-1)^{\theta(i\sqcup A',B')+d\cdot|B'|}
[c(i\sqcup A',x),c(B',y)]
+
\sum_{I'=A'\sqcup B'}
(-1)^{\theta(A',i\sqcup B')+d\cdot|i\sqcup B'|}
[c(A',x),c(i\sqcup B',y)]
\\=
\sum_{I=A\sqcup B}
(-1)^{\theta(A,B)+d\cdot|B|}
[c(A,x),c(B,y)].
\qedhere
\end{multline*}
\end{proof}
\begin{crl}
Let $I\subset[m],$ $I=I''\sqcup I',$ $I''<I'.$ Let $x,y\in\Gamma$ be homogeneous, and let $\mathcal{A}\subset 2^{I'}\times 2^{I'}$ be a family of pairs of subsets. Then
\begin{equation}
\label{eq:regrouping_c}
\sum_{\begin{smallmatrix}
I'=A'\sqcup B':\\
(A',B')\in\mathcal{A}
\end{smallmatrix}}
\!\!\!\!\!\!
(-1)^{\theta(A',B')+|B'|}\,
c(I'',[c(A',x),c(B',y)])
=
\!\!\!\!\!\!\!\!
\sum_{\begin{smallmatrix}
I=A\sqcup B:\\
(A\cap I',B\cap I')\in\mathcal{A}
\end{smallmatrix}}
\!\!\!\!\!\!\!\!
(-1)^{\theta(A,B)+|B|}\,
[c(A,x),c(B,y)].
\end{equation}
\end{crl}
\begin{proof}
It follows from \eqref{c(I,[x,y])_expanded} and from identities $c(A'',c(A',x))=c(A''\sqcup A',x),$ $\theta(A''\sqcup A',B''\sqcup B')=\theta(A'',B'')+\theta(A',B')+|A'|\cdot|B''|$ that are true for $A'',B''<A',B'.$
\end{proof}
\begin{prp}
Let $J\subset[m]$ and $i,j\in J$ such that $i<j$ and $J_{>j}\neq\varnothing.$ Then
\begin{multline}
\label{eq:commutator_hard_identity}
c(J\setminus ij,[u_i,u_j])=
(-1)^{|J_{>j}|}
c(J\setminus i,u_i)
-
(-1)^{|J_{>i}|}
c(J\setminus j,u_j)
\\+
\sum_{\begin{smallmatrix}
J\setminus ij = A\sqcup B:\\
A_{>i},B_{>j}\neq\varnothing
\end{smallmatrix}}
(-1)^{\theta(A,B)+|B|}
[c(A,u_i),c(B,u_j)].
\end{multline}
\end{prp}

\begin{proof}
Denote $P=J_{<j},$ $Q=J_{>i}\cap J_{<j},$ $R=J_{>j}.$ Hence $P<i<Q<j<R$ and $R\neq\varnothing.$ The left hand side is equal to $x:= c(P\sqcup Q,c(R,[u_i,u_j])).$ Denote also $y:=c(P\sqcup Q,[c(R,u_i),u_j)]),$ $z:=c(P\sqcup Q,[u_i,c(R,u_j)]).$ Then
\begin{align*}
x =&~
y
+(-1)^{|R|}z+
\sum_{\begin{smallmatrix}
R=A'\sqcup B':\\
A',B'\neq\varnothing
\end{smallmatrix}}
(-1)^{\theta(A',B')+|B'|}
c(P\sqcup Q,[c(A',u_i),c(B',u_j)])
\\
\underset{\eqref{eq:regrouping_c}}{=}&~
y+(-1)^{|R|}z+\sum_{\begin{smallmatrix}
J\setminus ij = A\sqcup B:\\
A_{>j},B_{>j}\neq\varnothing
\end{smallmatrix}}
(-1)^{\theta(A,B)+|B|}
[c(A,u_i),c(B,u_j)],
\end{align*}
$$y=(-1)^{|R|}c(P\sqcup Q,[u_j,c(R,u_i)])=(-1)^{|R|}c(J\setminus i,u_i),$$
\begin{multline*}
z=c(P,c(Q,[u_i,c(R,u_j)]))
\underset{\eqref{c(I,[x,y])_expanded}}{=}
c(P,[c(Q,u_i),c(R,u_j)])+(-1)^{|Q|}\underbrace{c(P,[u_i,c(Q\sqcup R,u_j)])}_{=c(J\setminus j,u_j)}
\\+
\sum_{\begin{smallmatrix}
Q=A_2\sqcup B_2:\\
A_2,B_2\neq\varnothing
\end{smallmatrix}}
(-1)^{\theta(A_2,B_2)+|B_2|}
c(P,[c(A_2,u_i),c(B_2\sqcup R,u_j)])
\\
\underset{\eqref{eq:regrouping_c}}{=}
(-1)^{|R|}\sum_{\begin{smallmatrix}
J\setminus ij = A\sqcup B:\\
B\cap Q=\varnothing,\\
A_{>j}=\varnothing\end{smallmatrix}
}
(-1)^{\theta(A,B)+|B|}
[c(A,u_i),c(B,u_j)]
\\+
(-1)^{|Q|}c(J\setminus j,u_j)
+
(-1)^{|R|}\sum_{\begin{smallmatrix}
J\setminus ij = A\sqcup B:\\
B\cap Q\neq\varnothing,Q;\\
A_{>j}=\varnothing
\end{smallmatrix}}
(-1)^{\theta(A,B)+|B|}
[c(A,u_i),c(B,u_j)].
\end{multline*}
Therefore,
\begin{multline*}
x=(-1)^{|R|}c(J\setminus i,u_i)+(-1)^{|Q|+|R|}c(J\setminus j,u_j)
\\+
\sum_{\begin{smallmatrix}
J\setminus ij = A\sqcup B:\\
A_{>i}\neq\varnothing,\\
A_{>j}=\varnothing\end{smallmatrix}
}(-1)^{
\theta(A,B)
+|B|
}[c(A,u_i),c(B,u_j)]+
\sum_{\begin{smallmatrix}
J\setminus ij=A\sqcup B:\\
A_{>j},B_{>j}\neq\varnothing
\end{smallmatrix}}
(-1)^{\theta(A,B)+|B|}[c(A,u_i),c(B,u_j)].
\end{multline*}
In the first sum the condition $B_{>j}\neq\varnothing$ is always true, since $R=A_{>j}\sqcup B_{>j},$ $A_{>j}=\varnothing$ and $R\neq\varnothing.$ In the second sum, $A_{>i}\neq\varnothing$ is always true. Hence the sums can be merged:
$$
\sum_{\begin{smallmatrix}
J\setminus ij = A\sqcup B:\\
A_{>i},B_{>j}\neq\varnothing\end{smallmatrix}
}(-1)^{
\theta(A,B)
+|B|
}[c(A,u_i),c(B,u_j)].
$$
Using $|R|=|J_{>j}|$ and $|Q+|R|=|J_{>i}|-1,$ we obtain \eqref{eq:commutator_hard_identity}.
\end{proof}

\end{document}